%% file: x-2-p-4-siu-yeung.tex
\documentclass[12pt,twoside,leqno,openany]{amsart}
\usepackage{amssymb,amsbsy,amsmath,amsfonts,amssymb,amscd,times,
graphics,color,xypic,footmisc,fancyhdr,multicol,fancybox,
graphicx,mathrsfs,rotating,ifthen,wasysym}
\usepackage[all]{xy}
\usepackage[french]{babel}
\usepackage[utf8]{inputenc}
\usepackage[T1]{fontenc}
\definecolor{black}{cmyk}{1.,1.,1.,1.0}
\definecolor{blue}{cmyk}{1.,1.,0.,0.63}
\definecolor{red}{cmyk}{0.,1.,1.,0.63}
\definecolor{green}{cmyk}{1.,0.,1.,0.63}

\input macros.tex

\let\mathcal\mathscr

\begin{document}

$\:$

\bigskip\bigskip

\begin{center}

{\large\bf Différentielles de jets de type Siu-Yeung}

\medskip

{\large\bf sur des surfaces $X^2$ intersections complètes 
dans $\P^4(\C)$}

\end{center}

\bigskip

\begin{center}
Jo\"el {\sc Merker}
\end{center}

\medskip

\begin{center}
\begin{minipage}[t]{11.75cm}
\baselineskip =0.35cm {\scriptsize

\centerline{\bf Table des matières}

\smallskip

{\bf 1.~Introduction
\dotfill~\pageref{introduction}.}

{\bf 2.~Paire de courbes algébriques complexes dans $\P^2(\C)$
\dotfill~\pageref{paire-courbes-P2}.}

{\bf 3.~Théorème d'injectivité
\dotfill~\pageref{theoreme-injectivite}.}

{\bf 4.~Géométrie des sections holomorphes de $\Sym^m\,T_X^*$
\dotfill~\pageref{sections-holomorphes}.}

{\bf 5.~Contraintes et degrés de liberté
\dotfill~\pageref{contraintes-liberte}.}

}\end{minipage}
\end{center}

\medskip

\hfill
\begin{minipage}[t]{10.25cm}
{\scriptsize\sf\em
The dominant theme is the {\em interplay} between the extrinsic projective
geometry of algebraic subvarieties of $\P^n(\C)$ and their intrinsic geometric 
features.}\hfill
Phillip~{\sc Griffiths}.
\end{minipage}


\bigskip

\centerline{\bf 1.~Introduction}
\label{introduction}
\HEAD{1.~Introduction}{
Jo\"el Merker, Département de Mathématiques d'Orsay}

\medskip

Dans sa thèse, Brotbek (\cite{Brotbek-2011-these,
Brotbek-2011-arxiv}) est parvenu à démontrer la conjecture de
Debarre (\cite{ Debarre-2005})
en dimension $n = 2$, et ce aussi, en {\em toute} codimension
$c \geqslant 2$, pour des surfaces {\em génériques} $X^2 \subset
\mathbb{ P}^{ 2 + c} (\C)$ intersections complètes dont les degrés
$d_1, \dots, d_c$ sont tous $\geqslant \frac{ 8c+18}{c-1}$.

En généralisant aux intersections complètes
quelconques une théorie intrinsèque développée auparavant 
pour les hypersurfaces, 
Brotbek a en outre établi
que si les degrés $d_1, \dots, d_c$ d'un système
de $c$ équations algébriques d'une intersection complète arbitraire
$X^n \subset \P^{ n+c} ( \C)$ sont tous supérieurs ou égaux à:
\[
\bigg[
2^{n-1}\,
\big(2n-2\big)\,
\frac{n^2}{n+c+1}\,
\binom{2n-1}{n}
+1
\bigg]
\binom{n}{{\sf Ent}(\frac{n}{2})}
\frac{(2n+c)!}{(n+c)!}\,
\frac{(c-n)!}{c!},
\]
alors ${\sf Sym}^m T_X^*$ possède des sections holomorphes globales
non identiquement nulles, sans plus, toutefois, 
d'information effective.  Pour ce
qui concerne ladite théorie intrinsèque, {\em confer} aussi Bérczi
(\cite{ Berczi-2012}), Darondeau (\cite{ Darondeau-2013}), Demailly
(\cite{ Demailly-1997, Demailly-2011}), Diverio (\cite{ Diverio-2009,
Diverio-Rousseau-2011}), Merker (\cite{ Merker-2010}), Mourougane
(\cite{ Mourougane-2012}), Pa\u{u}n (\cite{ Paun-2012}), Rousseau
(\cite{ Diverio-Rousseau-2011, Rousseau-2007}).

En s'inspirant de l'approche extrinsèque première de Siu-Yeung 1996 (\cite{
Siu-Yeung-1996}, {\em confer} aussi~\cite{ Siu-1995,
Siu-2002, Siu-2004, Siu-2012}), notre objectif
ici est de construire des sections holomorphes de $\Sym^m T_X^*$ plus
explicites en coordonnées. La plupart des techniques mises en {\oe}uvre se
généralisent sans grande difficulté à la dimension quelconque $n
\geqslant 2$, pour les intersections complètes de codimension aussi
quelconque $c \geqslant 2$, mais il est préférable, en vue
d'explorations ultérieures d'exemples amples, de se cantonner d'abord
à la dimension $2$ et à la codimension $2$, sans céder à la 
tentation\,\,---\,\,certes attrayante!\,\,---\,\,de la généralisation 
indicielle.

Soit donc $X^2 \subset \P^4 ( \C)$ une surface intersection complète
de degrés:
\[
1
\,\leqslant\, 
{\bf d}
\,\leqslant\, 
{\bf e}.
\]
Une application du Théorème~4 de Brückmann
dans~\cite{Bruckmann-1997} exprime 
explicitement la caractéristique d'Euler du fibré
des $m$-formes différentielles symétriques $\Sym^m\, T_X^*$:
\[
\!\!\!\!\!\!\!\!\!\!\!\!\!\!\!\!\!\!\!\!
\footnotesize
\aligned
\chi_{\sf Euler}
\big(
X^2,\,
\Sym^m\,T_X^*
\big)
=
\frac{1}{
1!\,2!\,3!\,4!}
&
\Big\{\,\,
m^3
\big[
24\,d^2\,e^2
-
120\,\big(d^2\,e+d\,e^2\big)
+
360\,e\,d
\big]
+
\\
&
+
m^2
\big[
-\,72\,d^2\,e^2
-
72\,\big(d^3\,e+d\,e^3\big)
+
360\,\big(d^2\,e+d\,e^2\big)
-
720\,d\,e
\big]
+
\\
&
+
m
\big[
-\,60\,d^2\,e^2
-
48\,\big(d^3\,e+d\,e^3\big)
+
300\,\big(d^2\,e+d\,e^2\big)
-
660\,d\,e
\big]
+
\\
&
+
36\,d^2\,e^2
+
24\,\big(d^3\,e+d\,e^3\big)
-
180\,\big(d^2\,e+d\,e^2\big)
+
420\,d\,e
\Big\}.
\endaligned
\]
Idéalement, la positivité de cette caractéristique devrait gouverner
l'existence de sections holomorphes globales de $\Sym^m T_X^*$, y
compris au moyen de certains calculs explicites en coordonnées qui
sont encore très peu explorés à ce jour, calculs qui fourniraient les
informations cruciales demeurant actuellement inaccessibles avec les
outils développées par la théorie intrinsèque depuis plus d'une
cinquantaine d'années en géométrie algébrique contemporaine,
limitation principielle qu'ont particulièrement bien mise en lumière
les travaux cités de Yum-Tong Siu.

\smallskip

Loin d'un tel horizon, nous nous contenterons de détailler
très soigneusement la démonstration d'un
résultat dans cet esprit.

\medskip\noindent{\bf Théorème.}
{\em
Sur une surface $X^2 \subset \P^4$ d'équations affines:}
\[
\aligned
z^d
&
=
R(x,y),
\\
t^e
&
=
S(x,y),
\endaligned
\]
{\em où $R \in \C [x, y]$ et $S \in \C [x, y]$ sont 
deux polynômes de degrés respectifs:}
\[
752
\,\leqslant\,
d
\,\leqslant\,
e
\,\leqslant\,
\frac{1}{648}\,
d^2
\] 
{\em qui satisfont un nombre fini de
dispositions géométriques génériques, pour tous polynômes:}
\[
A_{j,k,p,q}(x,y)
=
\sum_{h+i\leqslant a}\,
A_{j,k,p,q}^{h,i}\,
x^h\,y^i\,
\,\in\,\C[x,y]
\]
{\em de degrés:}
\[
\deg\,
A_{j,k,p,q}
\,\leqslant\,
a
\leqslant
d-4m,
\]
{\em la différentielle de jets méromorphe extrinsèque:}
\[
\frac{{\sf J}\big(x,y,x',y'\big)}{y^d\,z^{m(d-1)}\,t^{m(e-1)}}
\]
{\em où:}
\[
\aligned
{\sf J}\big(x,y,x',y'\big)
=
\sum_{j+k+p+q=m}\,
&
A_{j,k,p,q}(x,y)\,
\big(x'\big)^j\,\big(y'\big)^k\,
\\
&
\Big(
x'\,R_x(x,y)+y'\,R_y(x,y)
\Big)^p\,
\Big(
x'\,S_x(x,y)+y'\,S_y(x,y)
\Big)^q
\\
&
\big(R(x,y)\big)^{m-p}\,
\big(S(x,y)\big)^{m-q}
\endaligned
\]
{\em possède une restriction à la surface $X^2$:}
\[
\frac{{\sf J}(x,y,x',y')}{y^d\,z^{m(d-1)}\,t^{m(e-1)}}
\bigg\vert_{X^2}
\]
{\em qui est une section {\em holomorphe} du fibré des différentielles
symétriques intrinsèques:}
\[
\Sym^m T_X^*,
\]
{\em pourvu seulement que le numérateur-polynôme:}
\[
{\sf J}\big(x,y,x',y'\big)
\,\equiv\,
y^d\,\widetilde{\sf J}\big(x,y,x',y'\big)
\]
{\em soit divisible par $y^d$, ce qui impose aux
coefficients $A_{j,k,p,q}^{h,i} \in \C$
en nombre:}
\[
\frac{(a+1)(a+2)}{2}\,
\frac{(m+1)(m+2)(m+3)}{1\cdot 2\cdot 3}
\]
{\em de satisfaire un certain système d'équations linéaires dépendant
de $R$ et de $S$ dont l'espace des solutions est de dimension:}
\[
\dim\,H^0\big(X,\,\Sym^mT_X^*\big)
\,\geqslant\,
\frac{1}{93312}\,d^3
-
\frac{61}{7776}\,d^2
-
\frac{17}{108}\,d
-
\frac{28}{27},
\]
{\em lorsqu'on choisit:}
\[
m
:=
\Ent\,\frac{d}{12}.
\]

\medskip\noindent{\bf Remerciements.}
Une partie de ces constructions, notamment en dimension quelconque, a
été présentée au {\sl Hayama Symposium XV}, organisé par Katsutoshi
Yamanoi, Kengo Hirachi et Hajime Tsuji du 21 au 24 Juillet 2012.  Les
premiers éléments de ce texte liminaire ont pris forme lors du {\sl
Memorial Symposium: "Geometry and Analysis on Manifolds"}, en
souvenir du Professeur Shoshichi Kobayashi${}^\dag$, organisé par Takushiro
Ochiai, Keizo Hasegawa, Toshiki Mabuchi, Yoshiaki Maeda, Junjiro
Noguchi, Yoshihiko Suyama, Takashi Tsuboi, du 22 au 25 Mai 2013 dans
la {\sl Graduate School of Mathematical Sciences} de l'Université de
Tokyo. \`A cette occasion,
nous avons bénéficié d'échanges particulièrement intéressants
avec Junjiro Noguchi concernant des perspectives d'exemples amples.

\bigskip


\centerline{\bf 2.~Paire de courbes algébriques
complexes génériques dans $\P^2 (\C)$}
\label{paire-courbes-P2}
\HEAD{2.~Paire de courbes algébriques
complexes génériques dans $\P^2 (\C)$}{
Jo\"el Merker, Département de Mathématiques d'Orsay}

\medskip

\noindent{\bf Géométrie initiale.}
Sur $\P^4$, soient les coordonnées homogènes:
\[
\big[U\colon X\colon Y\colon Z\colon T\big],
\]
et, dans la carte affine initiale $\{ U \neq 0\}$, soient
les coordonnées affines
\[
x
:=
\frac{X}{U},
\ \ \ \ \ \ \ \ \ \ \ \ \ \ \ \ \
y
:=
\frac{Y}{U},
\ \ \ \ \ \ \ \ \ \ \ \ \ \ \ \ \
z
:=
\frac{Z}{U},
\ \ \ \ \ \ \ \ \ \ \ \ \ \ \ \ \
t
:=
\frac{T}{U}.
\]
Spécifiant comme <<\,horizontal\,>> le
\[
\C_{xy}^2
\,\subset\,
\C_{xyzt}^4,
\]
on se donne deux courbes algébriques complexes
d'équations polynomiales
\[
0
=
R(x,y)
\ \ \ \ \ \ \ \ \ \ \ \ \ \ \ \ \
\text{\rm et}
\ \ \ \ \ \ \ \ \ \ \ \ \ \ \ \ \
0
=
S(x,y),
\]
dont les degrés respectifs
\[
{\bf d}
:=
\deg\,R
\ \ \ \ \ \ \ \ \ \ \ \ \ \ \ \ \
\text{\rm et}
\ \ \ \ \ \ \ \ \ \ \ \ \ \ \ \ \
{\bf e}
:=
\deg\,S,
\]
satisfont, sans perte de généralité
\[
e\geqslant d\geqslant 1.
\]
Après l'action éventuelle d'un automorphisme affine de ce $\C_{x,y}^2$
horizontal, on peut supposer que
\[
\aligned
R(x,y)
&
=
\alpha\,x^d
+
\beta\,y^d
+
R_d^*
+
R_{d-1}
+\cdots+
R_0,
\\
S(x,y)
&
=
\gamma\,x^e
+
\delta\,y^e
+
S_e^*
+
S_{e-1}
+\cdots+
S_0,
\endaligned
\]
avec:
\[
\aligned
&
\alpha\neq 0,
\ \ \ \ \ \ \ \ \ \ \ \ \ 
\beta\neq 0,
\\
&
\gamma\neq 0,
\ \ \ \ \ \ \ \ \ \ \ \ \ 
\delta\neq 0,
\endaligned
\]
et:
\[
\aligned
&
R_d^*,\ \
R_{d-1},\ \
\dots,\ \
R_0,
\\
&
S_e^*,\ \
S_{e-1},\ \
\dots,\ \
S_0,
\endaligned
\]
homogènes en $(x, y)$ de degrés correspondants, $R_d^*$ et $S_e^*$
n'ayant pas de terme pur.

\medskip\noindent{\bf Hypothèse géométrique générique.}
{\em Les deux courbes sont lisses, et s'intersectent transversalement en 
exactement ${\bf d} \, {\bf e}$ points distincts de $\P_{xy}^2$
tous situés dans le $\C_{xy}^2$ affine.}

\medskip\noindent{\bf Proposition.}
{\em Alors la surface $X^2 \subset \P^4$ d'équations affines
\[
\aligned
z^d
&
=
R(x,y)
\\
t^e
&
=
S(x,y)
\endaligned
\]
est lisse, intersection complète dans $\P^4$.}

\proof
Tout d'abord, dans ce $\C_{xyzt}^4$, en un point où $R \neq 0 \neq S$,
elle consiste en $d\, e$ graphes lisses. Puis, en un point où $R = 0
\neq S$, la différentielle de $R$ est non nulle, donc $\{ z^d = R\}$
est lisse, et $e$ graphes s'empilent. Enfin, là où $R = 0 = S$, les
différentielles sont {\em à dessein} supposées indépendantes.

Pour terminer, le $(1/x)$- et le $(1/y)$-changements de carte
capturent tous les points de $X^2 \cap \P_{\infty,xyzt}^3$
et laissent inchangée la forme des équations.
\endproof


\bigskip

\centerline{\bf 3.~Théorème d'injectivité}
\label{theoreme-injectivite}
\HEAD{3.~Théorème d'injectivité}{
Jo\"el Merker, Département de Mathématiques d'Orsay}

\medskip

Inspiré par~\cite{Siu-Yeung-1996}, on introduit
\[
{\sf J}
:=
\sum_{j+k+p+q\,=\,{\bf m}
\atop
j,k,p,q\,\in\,\N}\,
A_{j,k,p,q}(x,y)\,
(x')^j\,(y')^k\,
\big(R'\big)^p\,(S'\big)^q\,
\big(R\big)^{{\bf m}-p}\,
\big(S\big)^{{\bf m}-q},
\]
l'entier $m \geqslant 1$ étant à choisir, les $A_{j,k,p,q} \in \C [x,y]$ 
étant de degrés:
\[
\deg\,
\big(A_{j,k,p,q}\big)
\leqslant
{\bf a}
\leqslant
{\bf d}-2,
\]
l'entier $a$ étant à choisir, $(x', y')$ étant les coordonnées 
cotangentes associées, d'où
\[
R'
=
x'\,R_x+y'\,R_y
\ \ \ \ \ \ \ \ \ \ \ \ \ \ \ \ \
\text{\rm et}
\ \ \ \ \ \ \ \ \ \ \ \ \ \ \ \ \
S'
=
x'\,R_x+y'\,R_y.
\]

\medskip\noindent{\bf Hypothèse géométrique générique.}
{\em Deux parmi six courbes 
\[
\aligned
&
\{R=0\},
\ \ \ \ \ \ \ \ \ \ \ \ \ \ \ \ \ 
\{R_x=0\},
\ \ \ \ \ \ \ \ \ \ \ \ \ \ \ \ \ 
\{R_y=0\},
\\
&
\{S=0\},
\ \ \ \ \ \ \ \ \ \ \ \ \ \ \ \ \ 
\{S_x=0\},
\ \ \ \ \ \ \ \ \ \ \ \ \ \ \ \ \ \
\{S_y=0\},
\endaligned
\]
s'intersectent transversalement dans $\P_{xy}^2$ en un 
nombre égal au produit de leurs degrés respectifs de points distincts
tous situés dans $\C_{xy}^2$, tandis qu'aucun triplet n'a
de point en commun dans $\P^2$.}

\medskip

On peut même supposer que chacune des ces six courbes intersecte aussi
transversalement la droite à l'infini. D'autres hypothèses
géométriques qui sont génériquement satisfaites devront
être admises pour garantir le bon déroulement des raisonnements.

\medskip\noindent{\bf Théorème d'injectivité.}
{\em 
Alors ${\sf J} \equiv 0$ dans $\C [x,y,x',y']$ si et seulement si tous
les $A_{j,k,p,q} \equiv 0$ dans $\C[x,y]$.}

\medskip

Admettons temporairement l'énoncé-outil suivant.

\medskip\noindent
{\bf Lemme d'annulation identique par restriction.}
{\em 
Si un polynôme $A \in \C[ x, y]$ de degré $\leqslant {\bf d} - 2$
s'annule en tous les points de l'une des neuf intersections
de dimension nulle:
\[
\aligned
0
&
=
A\big\vert_{R=S=0},
\ \ \ \ \ \ \ \ \ \ \ \ \ \ \ \ \ \ \ \ \ \ \ \ \
0
=
A\big\vert_{R_x=S=0},
\ \ \ \ \ \ \ \ \ \ \ \ \ \ \ \ \ \ \ \ \ \ \ \ \ \,
0
=
A\big\vert_{R_y=S=0},
\\
0
&
=
A\big\vert_{R=S_x=0},
\ \ \ \ \ \ \ \ \ \ \ \ \ \ \ \ \ \ \ \ \ \ \ \ 
0
=
A\big\vert_{R_x=S_x=0},
\ \ \ \ \ \ \ \ \ \ \ \ \ \ \ \ \ \ \ \ \ \ \ \ 
0
=
A\big\vert_{R_y=S_x=0},
\\
0
&
=
A\big\vert_{R=S_y=0},
\ \ \ \ \ \ \ \ \ \ \ \ \ \ \ \ \ \ \ \ \ \ \ \ 
0
=
A\big\vert_{R_x=S_y=0},
\ \ \ \ \ \ \ \ \ \ \ \ \ \ \ \ \ \ \ \ \ \ \ \ 
0
=
A\big\vert_{R_y=S_y=0},
\endaligned
\]
alors $A \equiv 0$ dans $\C[ x, y]$.}

\proof[Démonstration du théorème d'injectivité] Le développement de $0
\equiv {\sf J}$ utilise naturellement la formule du binôme:
\[
\!\!\!\!\!\!\!\!\!\!\!\!\!\!\!\!\!\!\!\!
\aligned
0
&
\equiv
\sum_{j+k+p+q\,=\,m}\,
A_{j,k,p,q}(x,y)\,
(x')^j\,(y')^k\,
\big(x'R_x+y'R_y\big)^p\,(x'S_x+y'S_y\big)^q\,
\big(R\big)^{m-p}\,
\big(S\big)^{m-q}
\\
&
\equiv
\sum_{j+k+p_1+p_2+q_1+q_2=m}\,
A_{j,k,p_1+p_2,q_1+q_2}\,
{\textstyle{\frac{(p_1+p_2)!}{p_1!\,p_2!}}}\,
{\textstyle{\frac{(q_1+q_2)!}{q_1!\,q_2!}}}\,
(x')^{j+p_1+q_1}\,(y')^{k+p_2+q_2}\,
\\
&
\ \ \ \ \ \ \ \ \ \ \ \ \ \ \ \ \ \ \ \ \ \ \ \ \ \ \ \ \ \ \ \ \ \ 
\ \ \ \ \ \ \ \ \ \ \ \ \ \ \ \ \ \ \ \ \ \ \ \ \ \ \ \ \ \ \ \ \ \ 
\ \ \ \ \ \ \ \ \ \ \ \ \ \ \ \ \
\big(R_x\big)^{p_1}\,\big(R_y\big)^{p_2}\,
\big(S_x\big)^{q_1}\,\big(S_y\big)^{q_2}\,
\\
&
\ \ \ \ \ \ \ \ \ \ \ \ \ \ \ \ \ \ \ \ \ \ \ \ \ \ \ \ \ \ \ \ \ \ 
\ \ \ \ \ \ \ \ \ \ \ \ \ \ \ \ \ \ \ \ \ \ \ \ \ \ \ \ \ \ \ \ \ \ 
\ \ \ \ \ \ \ \ \ \ \ \ \ \ \ \ \ \ 
\big(R\big)^{m-p_1-p_2}\,\big(S\big)^{m-q_1-q_2}.
\endaligned
\] 
Ensuite, il est nécessaire de réorganiser tout cela 
en monômes $(x')^\alpha\,
(y')^\beta$ avec $\alpha + \beta = m$, ce qui, en posant:
\[
\alpha
:=
j+p_1+q_1
\ \ \ \ \ \ \ \ \ \ \ \ \ \ \ \ \
\text{\rm et}
\ \ \ \ \ \ \ \ \ \ \ \ \ \ \ \ \
\beta
:=
k+p_2+q_2,
\]
donne l'annulation identique dans $\C [x, y, x', y']$:
\[
\!\!\!\!\!\!\!\!\!\!\!\!\!\!\!\!\!\!\!\!
\aligned
&
0
\equiv
\sum_{\alpha+\beta\,=\,m}\,
(x')^\alpha\,(y')^\beta\,
\sum_{j+p_1+q_1\,=\,\alpha}\,\,\,
\sum_{k+p_2+q_2\,=\,\beta}\,
{\textstyle{\frac{(p_1+p_2)!}{p_1!\,p_2!}}}\,
{\textstyle{\frac{(q_1+q_2)!}{q_1!\,q_2!}}}\,
A_{j,k,p_1+p_2,q_1+q_2}\,
\\
&
\ \ \ \ \ \ \ \ \ \ \ \ \ \ \ \ \ \ \ \ \ \ \ \ \ \ \ \ \ \ \ \ \ \ 
\ \ \ \ \ \ \ \ \ \ \ \ \ \ \ \ \ \ \ \ \ \ \ \ \ \ \ \ \ \ \ \ \ \ 
\ \ \ \ \ \ \ \ \ \ \ \ \ \ \ \ \ \ \ \ \ \ \ \ \ \ \ \ \ \ \
\big(R_x\big)^{p_1}\,\big(R_y\big)^{p_2}\,
\big(S_x\big)^{q_1}\,\big(S_y\big)^{q_2}\,
\\
&
\ \ \ \ \ \ \ \ \ \ \ \ \ \ \ \ \ \ \ \ \ \ \ \ \ \ \ \ \ \ \ \ \ \ 
\ \ \ \ \ \ \ \ \ \ \ \ \ \ \ \ \ \ \ \ \ \ \ \ \ \ \ \ \ \ \ \ \ \ 
\ \ \ \ \ \ \ \ \ \ \ \ \ \ \ \ \ \  \ \ \ \ \ \ \ \ \ \ \ \ \ \
\big(R\big)^{m-p_1-p_2}\,\big(S\big)^{m-q_1-q_2}.
\endaligned
\]
On en déduit donc, pour tous $\alpha + \beta = m$, les annulations
identiques dans $\C[x, y]$:
\[
\!\!\!\!\!\!\!\!\!\!\!\!\!\!\!\!\!\!\!\!
\boxed{\,
\aligned
&
0
\,\equiv\,
\sum_{j+p_1+q_1=\alpha}\,
\sum_{k+p_2+q_2=\beta}\,
{\textstyle{\frac{(p_1+p_2)!}{p_1!\,p_2!}}}\,
{\textstyle{\frac{(q_1+q_2)!}{q_1!\,q_2!}}}\,
A_{j,k,p_1+p_2,q_1+q_2}\,
\big(R_x\big)^{p_1}\,
\big(R_y\big)^{p_2}\,
\big(S_x\big)^{q_1}\,
\big(S_y\big)^{q_2}
\\
&
{\scriptstyle{(\forall\,\alpha\,+\,\beta\,=\,m)}}
\ \ \ \ \ \ \ \ \ \ \ \ \ \ \ \ \ \ \ \ \ \ \ \ \ \ \ \ \ \ \ \ \ \ 
\ \ \ \ \ \ \ \ \ \ \ \ \ \ \ \ \ \ \ \ \ \ \ \ \ \ \ \ \ \ \ \ \ \ 
\ \ \ \ \ \ 
\big(R\big)^{m-p_1-p_2}\,
\big(S\big)^{m-q_1-q_2}.
\endaligned\,}
\]

En particulier, pour $\beta = 0$, la deuxième somme disparaît:
\[
0
\equiv
\sum_{j+p_1+q_1\,=\,m}\,
A_{j,0,p_1,q_1}\
\big(R_x\big)^{p_1}\,\big(S_x\big)^{q_1}\,
\big(R\big)^{m-p_1}\,\big(S\big)^{m-q_1}.
\]
Dans une proposition ci-dessous,
on va montrer que ceci implique l'annulation de tous
ces $A_{j,0,p_1,q_1}$. 

Mais avant de synthétiser l'énoncé 
adéquat, supposons, par récurrence sur un certain entier
$k'$ avec $1 \leqslant k' \leqslant m$, que l'on sache:
\[
\aligned
&
0
\equiv
A_{j,0,p,q},
\ \ \ \ \ \ \ \ \ \ \ \ \
0
\equiv
A_{j,1,p,q},
\ \ \ \ \
\dots\dots\dots
\ \ \ \ \
0
\equiv
A_{j,k'-1,p,q},
\\
&
{\scriptstyle{(\forall\,j+p+q\,=\,m)}}
\ \ \ \ \ \ \ \ \ \ \ \ \
{\scriptstyle{(\forall\,j+p+q\,=\,m\,-\,1)}}
\ \ \ 
\dots\dots\dots
\ \ \ \ \
{\scriptstyle{(\forall\,j+p+q\,=\,m\,-\,(k'\,-\,1))}}.
\endaligned
\]
Une certaine {\sl triangularité sérendipitaire} va alors seconder
agréablement la récurrence.

En effet, si nous faisons $\beta = k'$ dans l'égalité encadrée:
\[
0
\equiv
\sum_{j+p_1+q_1\,=\,m-k'}\,\,
\sum_{k+p_2+q_2\,=\,k'}\,
\big(
\mathmotsf{mêmes termes}
\big)
\]
nous pouvons décomposer la deuxième somme en:
\[
=
\sum_{k+p_2+q_2\,=\,k'
\atop\ \ \ \ \ \ \ \ \ 
k\,=\,k'}\,
\big(\mathmotsf{un unique terme}\big)
+
\sum_{
k+p_2+q_2\,=\,k'
\atop\ \ \ \ \ \ \ \ \ \ \ \
k\,\leqslant\,k'-1}
\big(
\mathmotsf{mêmes termes}
\big),
\]
mais alors l'écriture complète des sommes et des termes:
\[
\aligned
0
&
\equiv
\sum_{j+p_1+q_1\,=\,m-k'}\,
A_{j,k',p_1,q_1}\,
\big(R_x\big)^{p_1}\,\big(S_x\big)^{q_1}\,
\big(R\big)^{m-p_1}\,\big(S\big)^{m-q_1}
+
\\
&
\ \ \ \ \
+
\sum_{j+p_1+q_1\,=\,m-k'}\,
\sum_{k+p_2+q_2\,=\,k'
\atop\ \ \ \ \ \ \ \ \ \ \ \  
k\,\leqslant\,k'-1}\,
{\textstyle{\frac{(p_1+p_2)!}{p_1!\,p_2!}}}\,
{\textstyle{\frac{(q_1+q_2)!}{q_1!\,q_2!}}}\,
\underbrace{A_{j,k,p_1+p_2,q_1+q_2}}_{
\text{\sf tous}\,\,\equiv\,0,\,\text{\sf par récurrence}}\,
\\
&
\ \ \ \ \ \ \ \ \ \ \ 
\big(R_x\big)^{p_1}\,\big(R_y\big)^{p_2}\,
\big(S_x\big)^{q_1}\,\big(S_y\big)^{q_2}\,
\big(R\big)^{m-p_1-p_2}\,\big(S\big)^{m-q_1-q_2}
\endaligned
\] 
fait voir que l'hypothèse de récurrence annule de manière
sérendipitaire le reste triangulaire, d'où découle
que l'équation identique
se simplifie alors en une équation 
qui possède la {\em même} forme:
\[
0
\equiv
\sum_{j+p_1+q_1\,=\,m-k'}\,
A_{j,k',p_1,q_1}\,
\big(R_x\big)^{p_1}\,\big(S_x\big)^{q_1}\,
\big(R\big)^{m-p_1}\,\big(S\big)^{m-q_1}
\]
que celle correspondant au choix initial $\beta = m$. 

Après une division par le polynôme non identiquement nul en facteur
commun $(R)^{ k'} \, (S)^{k'}$, et après la substitution notationnelle
$m - k' \longmapsto m$, le théorème d'injectivité se ramène donc bien
à un énoncé n'incorporant que $R_x$ et $S_x$, sans aucune apparition
de $R_y$ et $S_y$.
\endproof

\noindent{\bf Proposition.}
{\em Si des polynômes $A_{p, q} \in \C [x, y]$ de degrés: 
\[
\deg\,A_{p,q}
\,\leqslant\,
{\bf d}-1
\,\leqslant\,
{\bf e}-1
\]
satisfont:
\[
0
\,\equiv\,
\sum_{p+q\,\leqslant\,m}\,
A_{p,q}\,
\big(R_x\big)^p\,\big(S_x\big)^q\,
\big(R\big)^{m-p}\,\big(S\big)^{m-q},
\]
ils doivent être tous identiquement nuls.}

\medskip

Il importe de mentionner ici qu'une équation similaire dans
laquelle interviendraient aussi des puissances de $R_y$ et/ou
de $S_y$ n'impliquerait pas l'annulation des coefficients
correspondants.

\proof
\'Ecrivons l'hypothèse en spécifiant une extrémité de la sommation
triangulaire:
\[
\aligned
0
\equiv
A_{m,0}\,\big(R_x\big)^m\,\big(S\big)^m
+
\sum_{
\substack{
p+q\,\leqslant\,m
\\
\ \ \ \ \ \ \
p\,\leqslant\,m-1}}
A_{p,q}\,
\big(R_x\big)^p\,\big(S_x\big)^q\,
\big(R\big)^{m-p}\,\big(S\big)^{m-q}.
\endaligned
\]
Observant que le reste est multiple de $R$, restreignons
l'équation à la courbe $\{ R = 0\}$:
\[
0
\equiv
A_{m,0}\,
\big(R_x\big)^{m}\,\big(S\big)^{m}
\Big\vert_{\{R=0\}}.
\]
Or sur la courbe
lisse connexe $\{ R = 0\}$, les deux polynômes $R_x$ et $S$ 
ne sont pas identiquement nuls, donc on peut diviser:
\[
0
\equiv
A_{m,0}\,
\Big\vert_{\{R=0\}}.
\]
Comme $R$ est irréductible, et comme $\C [x, y]$ est un anneau
factoriel, $A_{ m, 0}$ doit être multiple de $R$, 
mais alors l'hypothèse sur le degré des
$A_{ p, q}$ a été faite pour qu'on
puisse en déduire:
\[
0
\equiv
A_{m,0}.
\]

Supposons, par induction sur un entier $p' \leqslant m - 1$,
que nous ayons déjà établi les annnulations identiques:
\[
0
\equiv
A_{p,q}
\ \ \ \ \ \ \ \ \ \ \ \ \
{\scriptstyle{(\forall\,p\,\geqslant\,p'\,+\,1,\,\,\,
\forall\,q\,=\,0\,\cdots\,m\,-\,p)}}.
\]
Dans la somme identiquement nulle par hypothèse, il ne reste
alors plus que des termes avec $p \leqslant p'$, somme
que l'on décompose alors en:
\[
\sum_{p+q\,\leqslant\,m
\atop
\ \ \ \,
p\,\leqslant\,p'}
\,=\,
\sum_{q=0}^{m-p'}\,
\bigg\vert_{p=p'}\,\,
+
\sum_{p+q\,\leqslant\,m
\atop
\ \ \ \ \ \ 
p\,\leqslant\,p'-1}.
\]
En extrayant la puissance de $R$ qui s'avère être
en facteur commun, ces deux s'écrivent très en détail:
\[
\!\!\!\!\!\!\!\!\!\!\!\!\!\!\!\!\!\!\!\!
\small
\aligned
0
&
\equiv
\zero{\big(R\big)^{m-p'}}\,
\bigg(
\sum_{q=0}^{m-p'}\,
A_{p',q}\,
\big(R_x\big)^{p'}\,\big(S_x\big)^{q}\,\big(S\big)^{m-q}
+
\sum_{p+q\,\leqslant\,m
\atop
\ \ \ \ \ \ 
p\,\leqslant\,p'-1}\,
\underbrace{
A_{p,q}\,\big(R_x\big)^p\,\big(S_x\big)^q\,
\big(R\big)^{p'-p}\,\big(S\big)^{m-q}}_{
\text{\sf multiple de}\,\,R\,=\,{\rm O}(R)}
\bigg).
\endaligned
\]
Bien entendu, puisque $R \not\equiv 0$, on peut diviser,
et ensuite, par restriction à $\{ R = 0\}$, on obtient:
\[
0
\equiv
\zero{\big(R_x\big)^{p'}}\,
\sum_{q=0}^{m-p'}\,
A_{p',q}\,\big(S_x\big)^q\,\big(S\big)^{m-q}
\Big\vert_{\{R=0\}}.
\]
\`A nouveau, on peut diviser par $R_x$
qui est non identiquement nul sur $\{ R = 0\}$, ce qui donne en expansion:
\[
\!\!\!\!\!\!\!\!\!\!\!\!\!\!\!\!\!\!\!\!
0
\equiv
\Big[
A_{p',0}\,\big(S\big)^m
+
A_{p',1}\,\big(S_x\big)^1\,\big(S\big)^{m-1}
+
A_{p',2}\,\big(S_x\big)^2\,\big(S\big)^{m-2}
+\cdots+
A_{p',m-p'}\,\big(S_x\big)^{m-p'}\,\big(S\big)^{p'}
\Big]
\bigg\vert_{\{R=0\}}.
\]
\`A présent, restreignons cette identité à l'ensemble
des $d ( e-1)$ points distincts $\{ R = S_x = 0\}$:
\[
0
=
A_{p',0}\,(S)^m
\big\vert_{\{R=S_x=0\}}.
\]
Puisque 
la disposition géométrique générique assure que $S$ ne s'annule
en aucun de ces $d(e-1)$ points, on déduit que $A_{ p', 0}$
doit s'y annuler, et le lemme d'annulation identique
par restriction délivre:
\[
0
\equiv
A_{p',0}.
\]
Il ne reste alors plus que: 
\[
\!\!\!\!\!\!\!\!\!\!\!\!\!\!\!\!\!\!\!\!
0
\equiv
\zero{\big(S_x\big)^1}
\Big[
A_{p',1}\,\big(S\big)^{m-1}
+
A_{p',2}\,\big(S_x\big)^1\,\big(S\big)^{m-2}
+\cdots+
A_{p',m-p'}\,\big(S_x\big)^{m-1-p'}\,\big(S\big)^{p'}
\Big]
\bigg\vert_{\{R=0\}}.
\]
On peut alors diviser par le $(S_x)^1$ en facteur commun, puis itérer
le raisonnement (exercice mental) pour conclure:
\[
0
\equiv
A_{p',1}
\equiv
A_{p',2}
\equiv\cdots\equiv
A_{p',m-p'},
\]
ce qui achève la démonstration par récurrence sur l'entier $p'$.
\endproof

\proof[Démonstration du lemme d'annulation identique par restriction]
En fait, les neuf cas en question sont chacun conséquence directe d'un 
énoncé unique. 
\endproof

\noindent{\bf Lemme.}
{\em \'Etant donné deux polynômes irréductibles $Q$ et $P$ de $\C[ x,
y]$ de degrés respectifs ${\bf q} \geqslant {\bf p} \geqslant 1$
dont les courbes projectives associées s'intersectent en exactement
$p\, q$ points de multiplicité $1$ tous situés dans le $\C_{xy}^2$
affine:
\[
\aligned
p\,q
&
=
\Card\,
\big\{P=0\big\}_{\C^2}
\cap
\big\{Q=0\big\}_{\C^2}
\\
&
=
\Card\,
\big\{P=0\big\}_{\P^2}
\cap
\big\{Q=0\big\}_{\P^2},
\endaligned
\]
si un polynôme $A \in \C[ x, y]$ de degré ${\bf a} \leqslant {\bf p} -
1$ s'annule en tous ces $p\,q$ points:
\[
0
=
A\big\vert_{\{P=Q=0\}_{\C^2}},
\]
alors $A \equiv 0$ dans $\C [x, y]$.}

\proof
Supposons par l'absurde que $A \not\equiv 0$.
Le polynôme $Q$ étant irréductible, la décomposition de $A$ en facteurs
premiers ne peut incorporer $Q$ pour des raisons de degré. Le théorème
classique de Bézout
(\cite{ Griffiths-1989}, p.~87), 
qui s'applique donc, assure alors, en tenant
compte des multiplicités, que:
\[
\aligned
a\,q
&
=
\Card\,
\big\{A=0\big\}_{\P^2}
\cap
\big\{Q=0\big\}_{\P^2}
\\
&
\geqslant
\Card\,
\big\{A=0\big\}_{\C^2}
\cap
\big\{Q=0\big\}_{\C^2}
\\
\explication{Hypothèse principale}
\ \ \ \ \ \ \ \ \ \ \ \ \ \ \
&
\geqslant
\Card\,
\big\{P=0\big\}_{\C^2}
\cap
\big\{Q=0\big\}_{\C^2}
\\
&
=
p\,q,
\endaligned
\]
ce qui, après division par $q$, apporte l'inégalité contradictoire $a
\geqslant p$, et conclut la preuve.
\endproof

\medskip\noindent{\bf Observations en vue de généralisations.}

\smallskip\noindent{\bf (1)}\,
On ne peut pas améliorer l'énoncé du lemme avec l'inégalité
${\bf a} \leqslant {\bf q} - 1$ en échangeant les rôles
de $Q$ et de $P$ dans ce raisonnement, car lorsque
${\bf p} \leqslant {\bf a} \leqslant {\bf q} - 1$, 
il se pourrait que $A \not\equiv 0$ soit divisible par
$P$, et l'application du théorème de Bézout s'effondre.

\smallskip\noindent{\bf (2)}\,
Dans le lemme, il suffit que $Q$ soit irréductible.

\smallskip\noindent{\bf (3)}\,
Il n'est pas nécessaire que tous les points
d'intersection entre les courbes projectives
soient situés dans le $\C_{ x, y}^2$ affine.

\smallskip\noindent{\bf (4)}\,
Dans l'application de
ce lemme à la proposition, on n'a en fait
besoin que du fait que:
\[
\big\{
R
=
0
=
S_x
\big\}
\]
consiste en ${\bf d} ( {\bf e} - 1)$ points distincts qui soient
situés dans $\C_{ x,y}^2$.


\bigskip

\centerline{\bf 4.~Géométrie des sections holomorphes de $\Sym^m\,T_X^*$}
\label{sections-holomorphes}
\HEAD{4.~Géométrie des sections holomorphes de $\Sym^m\,T_X^*$}{
Jo\"el Merker, Département de Mathématiques d'Orsay}

\medskip

Rappelons qu'après développement des puissances 
de $R' = x'\, R_x + y'\, R_y$ et de $S' = x'\, S_x + y'\, S_y$, 
le polynôme de jets ${\sf J}$ devient un certain
polynôme homogène de degré $m$ en $(x', y')$ 
dont les coefficients sont des formes
linéaires en les $A_{j,k,p,q}$:
\[
{\sf J}
=
\sum_{\alpha+\beta\,=\,m}\,
(x')^\alpha\,(y')^\beta\,
\Lambda_{\alpha,\beta}
\big(A_\bullet;\,
R_x,R_y,S_x,S_y,R,S
\big),
\]
où:
\[
\label{deja-vu}
\boxed{\,
\aligned
\Lambda_{\alpha,\beta}
&
:=
\sum_{j+p_1+q_1\,=\,\alpha}\,\,\,
\sum_{k+p_2+q_2\,=\,\beta}\,
{\textstyle{\frac{(p_1+p_2)!}{p_1!\,p_2!}}}\,
{\textstyle{\frac{(q_1+q_2)!}{q_1!\,q_2!}}}\,
A_{j,k,p_1+p_2,q_1+q_2}\,
\\
&
\ \ \ \ \ \ \ \ \ \ \ \ \ \ \ \ \ \ \ \ \ \ \ \ \ \ \ \ \ \ \ \ \ \ 
\ \ \ \ \ \ \ \ \ \ 
\big(R_x\big)^{p_1}\,\big(R_y\big)^{p_2}\,
\big(S_x\big)^{q_1}\,\big(S_y\big)^{q_2}\,
\\
&
\ \ \ \ \ \ \ \ \ \ \ \ \ \ \ \ \ \ \ \ \ \ \ \ \ \ \ \ \ \ \ \ \ \
\ \ \ \ \ \ \ \ \ \  
\big(R\big)^{m-p_1-p_2}\,\big(S\big)^{m-q_1-q_2},
\endaligned}
\]
coefficients que nous abrégerons parfois simplement en:
\[
\Lambda_{\alpha,\beta}(x,y).
\]

\'Etant donné un entier $c \geqslant 1$ dont la valeur exacte:
\[
{\bf c}
\,\approx\,
{\bf d} 
\]
pourrait varier en restant
proche de $d$, on introduit le quotient suivant:
\[
\aligned
\frac{{\sf J}}{y^c\,z^{m(d-1)}\,t^{m(e-1)}}
=
\frac{1}{z^{m(d-1)}\,t^{m(e-1)}}
\sum_{\alpha+\beta\,=\,m}\,
(x')^\alpha\,(y')^\beta\,
\frac{\Lambda_{\alpha,\beta}(x,y)}{y^c},
\endaligned
\]
qui est une différentielle de jets {\em méromorphe et
extrinsèque} définie sur $\C^4$ tout entier.

La Section~5 ci-dessous sera consacrée à établir, 
par de simples arguments de comptage, que pour toute paire
de polynômes $R$ et $S$ en position général de degrés respectifs:
\[
{\bf d}
\,\leqslant\,
{\bf e}
\,\leqslant\,
{\sf constante}\cdot{\bf d}^2,
\]
il existe ${\bf m}$ assez grand pour que le {\sl système linéaire
de divisibilité} (contraintes):
\[
\mathmotsf{tous les}\,\,\Lambda_{\alpha,\beta}(x,y)
\,\equiv\,
y^c\,
\Lambda_{\alpha,\beta}^\sim(x,y)
\]
d'inconnues les {\em coefficients} des polynômes
(degrés de liberté):
\[
A_{j,k,p,q}(x,y)
=
\sum_{h+i\,\leqslant\,{\bf a}}\,
\underbrace{A_{j,k,p,q}^{i,h}}_{\in\,\C}\,
x^h\,y^i
\]
possède de nombreuses solutions.

Dans ces conditions, lorsque les $A_\bullet$ appartiennent à cet
espace de solutions, la différentielle de jets méromorphe
extrinsèque:
\[
\aligned
\frac{{\sf J}}{y^c\,z^{m(d-1)}\,t^{m(e-1)}}
=
\frac{1}{z^{m(d-1)}\,t^{m(e-1)}}
\sum_{\alpha+\beta\,=\,m}\,
(x')^\alpha\,(y')^\beta\,
\Lambda_{\alpha,\beta}^\sim(x,y),
\endaligned
\]
s'avère être {\em holomorphe} dans:
\[
\C^4
\cap
\{z\neq 0\}
\cap
\{t\neq 0\},
\]
donc holomorphe en restriction à la surface dans ce même ouvert:
\[
\boxed{\,
\frac{{\sf J}}{y^c\,z^{m(d-1)}\,t^{m(e-1)}}
\bigg\vert_{
X^2
\cap 
\C^4
\cap
\{z\neq 0\}
\cap
\{t\neq 0\}}
\,\,\text{\rm holomorphe}.\,}
\]

\`A la fin de la présente 
Section 4, on verra aussi qu'un choix approprié de
${\bf m}$, ${\bf d}$, ${\bf e}$ assure que cette différentielle de
jets extrinsèque s'annule identiquement sur l'hyperplan à l'infini
$\P_\infty^3$, donc est holomorphe en restriction à $X^2 \cap
\P_\infty^3$.

\medskip\noindent{\bf Première proposition principale.}
{\em
Lorsque les $A_{j,k,p,q} (x,y)$ sont soumis aux contraintes de
divisibilité de tous les $\Lambda_{\alpha, \beta} ( x,y)$ par $y^c$,
la différentielle symétrique de jets méromorphe extrinsèque définie
sur $\C^4${\em :}
\[
\frac{{\sf J}}{y^c\,z^{m(d-1)}\,t^{m(e-1)}}
\]
s'avère possèder une restriction à la surface affine $X^2
\cap \C^4$ d'équations:
\[
\aligned
z^d
&
=
R(x,y),
\\
t^e
&
=
S(x,y),
\endaligned
\]
qui devient {\em holomorphe} comme différentielle symétrique de $1$-jets 
intrinsèque à $X^2 \cap \C^4$.}

\proof
On vient de voir qu'en dehors des deux hyperplans $\{ z = 0\}$
et $\{ t = 0 \}$, sa restriction à $X^2 \cap \C^4$ est holomorphe.
Il reste à étudier l'holomorphie
en un point quelconque:
\[
(x_0,y_0,z_0,t_0)
\] 
appartenant à l'un des trois ouverts d'espace linéaire:

\medskip\noindent{\bf Cas 1:}\,
$\{ z = 0 \} \cap \{ t \neq 0 \}$, qui intersecte $X^2 \cap \P^4$
en une courbe;

\medskip\noindent{\bf Cas 2:}\,
$\{ z \neq 0 \} \cap \{ t = 0 \}$, qui intersecte $X^2 \cap \P^4$
aussi en une courbe;

\medskip\noindent{\bf Cas 3:}\,
$\{ z = 0 \} \cap \{ t = 0 \}$, qui intersecte $X^2 \cap \P^4$
en $d\, e$ points distincts.

\medskip

Localisons donc pour commencer l'étude de l'holomorphie
au voisinage d'un point appartenant au premier ensemble:
\[
(x_0,y_0,0,t_0)
\ \ \ \ \ \ \ \ \ \ \ \ \
{\scriptstyle{(t_0\,\neq\,0)}}.
\] 
Clairement, l'intersection $X^2 \cap \C_{xyt}^3
\cap \{ t \neq 0\}$ est la courbe lisse d'équations:
\[
\aligned
0
&
=
R(x,y),
\\
t^d
&
=
S(x,y),
\endaligned
\]
et puisque ${\sf J}$ est divisée par $z^{ m(d - 1)}$, le quotient
semble singulier le long de cette courbe. Or il n'en est rien! 

En effet, si on différentie une fois la première équation de $X^2$:
\[
z'\,d\,z^{d-1}
=
R',
\]
et si l'on remplace les occurence de la lettre $R$ dans 
${\sf J}$\,\,---\,\,ce qui correspond à effectuer un changement de carte
dans le fibré des jets intrinsèques à $X^2$\,\,---,
on constate que la singularité en $z$ disparaît:
\[
\!\!\!\!\!\!\!\!\!\!\!\!\!\!\!\!\!\!\!\!
\footnotesize
\aligned
\frac{{\sf J}}{y^c\,z^{m(d-1)}\,t^{m(e-1)}}
&
=
\frac{1}{y^c}\,
\frac{1}{t^{m(e-1)}}\,
\sum_{j+k+p+q\,=\,m}\,
A_{j,k,p,q}\,
(x')^j\,(y')^k\,
\frac{\big(z'\,d\,z^{d-1}\big)^p\,\big(z^d\big)^{m-p}}{z^{m(d-1)}}\,
\big(S'\big)^q\,\big(S\big)^{m-q}
\\
&
=
\frac{1}{y^c}\,
\underbrace{\frac{1}{t^{m(e-1)}}}_{
\text{\sf holomorphe}
\atop
\text{\sf puisque}\,\,t_0\,\neq\,0}\,
\sum_{j+k+p+q\,=\,m}\,
d^p\,A_{j,k,p,q}\,
(x')^j\,(y')^k\,(z')^p\,
\underbrace{\,\,\,z^{m-p}\,\,\,}_{
\text{\sf disparition}
\atop
\text{\sf de la singularité}}\,
\big(S'\big)^q\,\big(S\big)^{m-q}.
\endaligned
\]
Mais puisque ${\sf J}$ est aussi divisée par $y^c$, deux sous-cas sont
à considérer séparément:

\medskip\noindent$\text{\bf Sous-cas 1}_{\text{\bf 1}}$:\,
$y_0\neq 0$;

\medskip\noindent$\text{\bf Sous-cas 1}_{\text{\bf 2}}$:\,
$y_0 = 0$.

\medskip

Dans le premier sous-cas, $\frac{ 1}{ y^c}$ est holomorphe au 
voisinage de $y_0$, donc l'expression complète est manifestement
holomorphe. Toutefois, l'écriture de cette expression n'est pas encore 
finalisée {\em intrinsèquement}, puisqu'elle incorpore
les quatre coordonnées extrinsèques $(x, y, z, t)$, alors
que $X^2$ est $2$-dimensionnelle. Il est avisé, alors, d'achever
l'expression explicite de la restriction à $X^2$.

Au voisinage d'un tel point:
\[
(x_0,y_0,0,t_0)
\,\in\,X^2\cap\C^4
\ \ \ \ \ \ \ \ \ \ \ \ \
{\scriptstyle{(y_0\,\neq\,0,\,\,\,t_0\,\neq\,0)}},
\]
la surface doit être représentée comme
graphe d'une application holomorphe
locale $\C^2 \longrightarrow \C^2$.
La matrice jacobienne de ses deux équations polynomiales:
\[
\left(\!\!
\begin{array}{cccc}
-R_x & -R_y & dz^{d-1} & 0
\\
-S_x & -S_y & 0 & et^{e-1}
\end{array}
\!\!\right)
\]
a pour valeur en un tel point où $z_0 = 0$:
\[
\left(\!\!
\begin{array}{cccc}
-R_x^0 & -R_y^0 & 0 & 0
\\
-S_x^0 & -S_y^0 & 0 & e\,t_0^{e-1}
\end{array}
\!\!\right),
\]
et donc deux éventualités sont à considérer:
\[
R_x(x_0,y_0)
\,\neq\,0
\ \ \ \ \ \ \ \ \ \ \ \ \ \ \ \ \
\text{\rm ou}
\ \ \ \ \ \ \ \ \ \ \ \ \ \ \ \ \
R_y(x_0,y_0)
\,\neq\,0.
\]
(Géométriquement parlant, puisque $(z_0)^d = 0$, 
le point $(x_0, y_0)$ appartient
à la courbe lisse $\{ R = 0\}$ vue dans le $\C_{xy}^2$
qui est contenu dans le $\C_{xyt}^3$ d'intersection
$\C^4 \cap \{ z = 0\}$ dont on est parti pour 
étudier le Cas 1.) 

Pour fixer les idées, supposons que $R_y^0 \neq 0$, l'autre
éventualité étant symétrique.
Le théorème des fonctions implicites fournit alors deux
applications holomorphes locales:
\[
{\tt Y}
\colon\ \ \ 
\C^2
\longrightarrow
\C
\ \ \ \ \ \ \ \ \ \ \ \ \ \ \ \ \
\text{\rm et}
\ \ \ \ \ \ \ \ \ \ \ \ \ \ \ \ \
{\tt T}
\colon\ \ \ 
\C^2
\longrightarrow
\C
\] 
définies dans un certain voisinage de $(x_0, 0)$ dans
$\C_{xz}^2$ satisfaisant:
\[
{\tt Y}(x_0,z_0)
=
y_0
\ \ \ \ \
{\scriptstyle{(\neq\,0)}}
\ \ \ \ \ \ \ \ \ \ \ \ \ \ \ \ \
\text{\rm et}
\ \ \ \ \ \ \ \ \ \ \ \ \ \ \ \ \
{\tt T}(x_0,z_0)
=
t_0
\ \ \ \ \
{\scriptstyle{(\neq\,0)}}
\]
qui représentent $X^2$ comme le graphe:
\[
\aligned
y
&
=
{\tt Y}(x,z),
\\
t
&
=
{\tt T}(x,z).
\endaligned
\]
Autrement dit, les deux équations:
\[
\aligned
z^d
&
\equiv
R\big(x,{\tt Y}(x,z)\big),
\\
\big(T(x,z)\big)^e
&
\equiv
S\big(x,{\tt Y}(x,z)\big),
\endaligned
\]
sont satisfaites identiquement pour tout $(x,z)$
dans le voisinage de $(x_0,z_0)$ en question. 

\smallskip

Ensuite, la différentiation directe 
de ces deux équations graphées donne: 
\[
\aligned
y'
&
=
x'\,{\tt Y}_x
+
z'\,{\tt Y}_z,
\\
t'
&
=
x'\,{\tt T}_x
+
z'\,{\tt T}_z.
\endaligned
\]

Maintenant, en revenant à la différentielle symétrique
méromorphe extrinsèque, son {\em expression induite
achevée} est naturellement celle que l'on obtient
en effectuant tous les remplacements nécessaires
pour obtenir une expression n'incorporant que
$(x, z, x', z')$:
\[
\footnotesize
\aligned
\frac{{\sf J}}{y^c\,z^{m(d-1)}\,t^{m(e-1)}}
&
=
\underbrace{\frac{1}{{\tt Y}(x,z)^{c}}}_{
\text{\sf localement holomorphe}
\atop
\text{\sf puisque}\,\,{\tt Y}(x_0,z_0)\,\neq\,0}\,
\underbrace{\frac{1}{{\tt T}(x,z)^{m(e-1)}}}_{
\text{\sf localement holomorphe}
\atop
\text{\sf puisque}\,\,{\tt T}(x_0,z_0)\,\neq\,0}\,
\sum_{j+k+p+q\,=\,m}\,
d^p\,
A_{j,k,p,q}\big(x,{\tt Y}(x,z)\big)\,
\\
&
\ \ \ \ \ \
(x')^j\,
\big(x'\,{\tt Y}_x(x,z)+z'\,{\tt Y}_z(x,z)\big)^k\,(z')^p\,
\underbrace{\,\,\,z^{m-p}\,\,\,}_{
\text{\sf singularité}
\atop
\text{\sf éliminée}}\,
\\
&
\ \ \ \ \ \
\Big(
x'\,S_x\big(x,{\tt Y}(x,z)\big)
+
\big(x'\,{\tt Y}_x(x,z)+z'\,{\tt Y}_z(x,z)\big)\,
S_y\big(x,{\tt Y}(x,z)\big)
\Big)^q\,
\\
&
\ \ \ \ \ \
\Big(S\big(x,{\tt Y}(x,z)\big)\Big)^{m-q}.
\endaligned
\]
{\em On voit bien alors que cette expression est effectivement
holomorphe dans un voisinage du $2$-plan-fibre}:
\[
\{x_0\}\times\{y_0\}\times\C_{x'}\times\C_{z'}
\,\subset\,
\C_{xzx'z'}^4.
\]
Bien entendu, cette expression pourrait encore être développée
pour la réorganiser en polynôme homogène de degré $m$ en
les deux variables de jets $(x', z')$.

\bigskip

Ensuite, pour ce qui est du deuxième
$\text{\rm Sous-cas 1}_{\text{\rm 2}}$, 
soit donc un point quelconque:
\[
(x_0,0,0,t_0)
\,\in\,
X^2
\cap 
\C^4
\ \ \ \ \ \ \ \ \ \ \ \ \
{\scriptstyle{(t_0\,\neq\,0)}}.
\]
En un tel point, la matrice Jacobienne des deux équations
affines de $X^2$ a la même expression:
\[
\left(\!\!
\begin{array}{cccc}
-R_x^0 & -R_y^0 & 0 & 0
\\
-S_x^0 & -S_y^0 & 0 & e\,t_0^{e-1}
\end{array}
\!\!\right).
\]

\medskip\noindent{\bf Disposition géométrique générique.}
{\em
La droite $\{ y = 0\}$ dans le plan affine $\C_{xy}^2 \subset
\P_{xy}^2$ intersecte les deux courbes $\{ R = 0\}$ et $\{ S = 0\}$
en, respectivement, $d$ et $e$ points de multiplicité $1$ où
les {\em deux} dérivées partielles $R_x$, $R_y$ et $S_x$, $S_y$ ne
s'annulent pas.}

\medskip

On peut donc à nouveau appliquer le théorème des fonctions implicites
et représenter de la même façon $X^2$, au moyen de fonctions
graphantes ${\tt Y}$ et ${\tt T}$ qui sont {\em a priori} distinctes
des précédentes sans toutefois introduire de nouvelle notation
différente, localement au voisinage d'un tel point $(x_0, 0, 0, t_0)$
comme:
\[
\aligned
y
&
=
{\tt Y}(x,z),
\\
t
&
=
{\tt T}(x,z).
\endaligned
\]
Bien entendu, ${\tt Y}(x_0,0) = 0$ et ${\tt T} ( x_0, 0) = t_0 \neq 
0$.

Ensuite, nous affirmons que l'intersection locale entre les deux
courbes $\{ z = 0\}$ et $\{ {\tt Y} ( x, z) = 0\}$ se réduit au point
$(x_0, 0)$, et qu'en fait:
\[
{\tt Y}_x(x_0,0)
\neq
0,
\]
ce qui assure alors un contrôle de la multiplicité,
égale à $1$, cette intersection.

\smallskip

En effet, les deux fonctions graphantes ${\tt Y}$ et ${\tt T}$ satisfont 
bien entendu les deux identités:
\[
\aligned
z^d
&
\equiv
R\big(x,{\tt Y}(x,z)\big),
\\
\big(T(x,z)\big)^e
&
\equiv
S\big(x,{\tt Y}(x,z)\big).
\endaligned
\]
Si l'on différentie la première par rapport à $x$, on obtient
classiquement l'expression de la dérivée partielle:
\[
{\tt Y}_x(x,z)
=
-\,
\frac{R_x\big(x,{\tt Y}(x,z)\big)}{
R_y\big(x,{\tt Y}(x,z)\big)},
\]
ce qui, au point $(x_0, 0)$, donne:
\[
\aligned
{\tt Y}_x(x_0,0)
=
-\,
\frac{R_x\big(x_0,0\big)}{
R_y\big(x_0,0\big)}
\neq
0,
\endaligned
\]
quantité non nulle grâce à la disposition géométrique générique
qui a été effectuée à l'avance.

\medskip\noindent{\bf Assertion.}
{\em 
Au voisinage de $(x_0, 0, 0, t_0)$, avec $t_0 \neq 0$, la restriction
à $X^2$ de la différentielle symétrique méromorphe définie sur $\C^4$:
\[
\frac{{\sf J}(x,y,x',y')}{
y^c\,z^{m(d-1)}\,t^{m(e-1)}}
=
\frac{{\sf J}(x,y,x',y')}{
{\tt Y}(x,z)^c\,z^{m(d-1)}\,{\tt T}(x,z)^{m(e-1)}}
\]
envisagée dans les coordonnées intrinsèques locales $(x, z, x', z')$
sur $X^2$ et qui est manifestement holomorphe en dehors de la réunion
des lieux polaires:
\[
\big\{(x,z,x',z')\colon\,
(x,z)\,\,\text{\rm près de}\,\,(x_0,0),\,\,\,
x'\in\C,\,z'\in\C
\big\}
\Big\backslash
\Big(
\big\{{\tt Y}(x,z)=0\big\}
\cup
\{z=0\}
\Big)
\] 
se prolonge en fait holomorphiquement à travers chacune des deux
courbes épointées de $\C_{xz}^2${\em :}
\[
\Big(
\big\{{\tt Y}(x,z)=0\big\}
\big\backslash\{z=0\}
\Big)
\bigcup
\Big(
\{z=0\}
\big\backslash
\big\{{\tt Y}(x,z)=0\big\}
\Big).
\]
De plus, l'intersection au voisinage de $(x_0, 0)$ 
entre ces deux courbes se réduit à{\em :}
\[
\big\{Y(x,z)=0\big\}
\cap
\{z=0\}
=
\{(x_0,0)\},
\]
qui est un seul point de multiplicité $1$, 
et le théorème de Riemann-Hurwitz-Hartogs
d'élimination des singularités 
pour les fonctions holomorphes locales définies en dehors
d'une sous-variété complexe de codimension $\geqslant 2$
assure qu'au voisinage d'un tel point:
\[
(x_0,0,0,y_0)
\ \ \ \ \ \ \ \ \ \ \ \ \
{\scriptstyle{(t_0\,\neq\,0)}},
\] 
la restriction à $X^2$ de ladite différentielle symétrique
méromorphe sur $\C^4$ est holomorphe, dans les
coordonnées intrinsèques locales $(x, z, x', z')$, sur
tout un voisinage de ce point-base `croix' la fibre:}
\[
\{(x_0,0)\}\times\C_{x'z'}^2.
\]

\proof
Rappelons que les coefficients $A_{ j, k, p, q} ( x, y)$ vont
être choisis dans la 
Section suivante pour qu'une propriété de divisibilité par $y^c$ de
$J ( x, y, x', y')$ soit satisfaite:
\[
J(x,y,x',y')
\equiv
y^c\,\widetilde{J}(x,y,x',y'),
\] 
une telle identité valant entre polynômes de $\C[ x, y, x', y']$.
Maintenant, effectuons la restriction intrinsèque nécessaire, en
observant que le $y^c$ singulier au dénominateur est alors effacé
par le $y^c$ au numérateur:
\[
\footnotesize
\aligned
\frac{{\sf J}(x,y,x',y')}{
y^c\,z^{m(d-1)}\,t^{m(e-1)}}
\bigg\vert_{X^2\,\,\text{\rm près de}\,\,(x_0,0,0,t_0)}
&
=
\frac{\zero{y^c}\,\widetilde{\sf J}(x,y,x',y')}{
\zero{y^c}\,z^{m(d-1)}\,t^{m(e-1)}}
\bigg\vert_{X^2\,\,\text{\rm près de}\,\,(x_0,0,0,t_0)}
\\
&
=
\frac{\widetilde{\sf J}
\big(
x,\,{\tt Y}(x,z),\,x',\,x'\,{\tt Y}_x(x,z)+z'\,{\tt Y}_z(x,z)\big)}{
z^{m(d-1)}\,{\tt T}(x,z)^{m(e-1)}}
\bigg\vert_{(x,z)\,\,\text{\rm près de}\,\,(x_0,0)}.
\endaligned
\]
Puisque ${\tt T}(x_0,0) = t_0 \neq 0$, et puisqu'on
divise par $z^{ m ( d-1)}$
cette expression est holomorphe au voisinage de
tout point de $X^2$ de coordonnées
intrinsèques $(x_1, z_1)$ avec $z_1 \neq 0$.
En particulier, elle est holomorphe en tout point:
\[
(x_1,z_1)
\,\in\,
\big\{{\tt Y}(x,z)=0\big\}
\big\backslash
\{z=0\}.
\]

\smallskip

Ensuite, le même calcul de restriction intrinsèque que 
celui que nous avons effectué
il y a quelques instants au voisinage d'un point $(x_0, y_0, 0, t_0)$,
donne\,\,---\,\,avec d'autres fonctions graphantes ${\tt Y}$ et 
${\tt T}$\,\,---:
\[
\footnotesize
\aligned
\frac{{\sf J}(x,z,x',z')}{y^c\,z^{m(d-1)}\,t^{m(e-1)}}
\bigg\vert_{X^2\,\,\text{\rm près de}\,\,(x_0,0,0,t_0)}
&
=
\frac{1}{{\tt Y}(x,z)^{c}}\,
\underbrace{\frac{1}{{\tt T}(x,z)^{m(e-1)}}}_{
\text{\sf localement holomorphe}
\atop
\text{\sf puisque}\,\,{\tt T}(x_0,0)\,=\,t_0\,\neq\,0}\,
\sum_{j+k+p+q\,=\,m}\,
d^p\,
A_{j,k,p,q}\big(x,{\tt Y}(x,z)\big)\,
\\
&
\ \ \ \ \ \
(x')^j\,
\big(x'\,{\tt Y}_x(x,z)+z'\,{\tt Y}_z(x,z)\big)^k\,(z')^p\,
\underbrace{\,\,\,z^{m-p}\,\,\,}_{
\text{\sf singularité}
\atop
\text{\sf éliminée}}\,
\\
&
\ \ \ \ \ \
\Big(
x'\,S_x\big(x,{\tt Y}(x,z)\big)
+
\big(x'\,{\tt Y}_x(x,z)+z'\,{\tt Y}_z(x,z)\big)\,
S_y\big(x,{\tt Y}(x,z)\big)
\Big)^q\,
\\
&
\ \ \ \ \ \
\Big(S\big(x,{\tt Y}(x,z)\big)\Big)^{m-q}.
\endaligned
\]
Visiblement, le lieu singulier (rémanent) du membre de droite est 
précisément la courbe de $\C_{xz}^2$ d'équation:
\[
\big\{
{\tt Y}(x,z)=0
\big\},
\]
qui passe par le point $(x_0, 0)$.  En particulier, cette deuxième
expression est holomorphe en tout point:
\[
(x_1,z_1)
\,\in\,
\{z=0\}
\big\backslash
\big\{{\tt Y}(x,z)=0\big\}
\]
proche de $(x_0, 0)$.

\smallskip

Ainsi, lorsqu'on effectue la {\em synthèse entre ces deux expressions
intrinsèques complémentaires de la \underline{même} différentielle
symétrique méromorphe extrinsèque}, on constate que ses
singularités apparentes s'évanouissent à travers tout point qui n'est
pas situé sur l'intersection entre les deux courbes de $\C_{ xz}^2$:
\[
\big\{{\tt Y}(x,z)=0\big\}
\cap
\{y=0\},
\]
et donc ses singularités rémanentes sont confinées à être contenues
dans:
\[
\big\{z={\tt Y}(x,z)=0\}
\times
\C_{x'y'}^2,
\]
et comme avons déjà fait
observer à l'avance que l'intersection locale entre ces deux
courbes $\{ z = 0\}$ et $\{ {\tt Y} ( x, z) = 0\}$ se réduit au point
$(x_0, 0)$, la démonstration de l'Assertion se termine 
grâce à une application du théorème de
Riemann-Hurwitz-Hartogs d'élimination des singularités
de codimension $\geqslant 2$.
\endproof

\noindent{\bf Résumé intermédiaire.}
Ainsi dans le Cas~1, l'holomorphie de ${\sf J} \big/ y^c \, z^{
m(d-1)} \, t^{ m ( e - 1)}$ est complètement établie au
voisinage de tout point $(x_0, y_0, 0, t_0)$ appartenant à $\{ z = 0
\} \cap \{ t \neq 0\}$.

\medskip

Le traitement du Cas~2 s'exécute d'une manière parfaitement
symétrique, puisque les deux équations de $X^2$ sont de la même forme,
et on obtient que la différentielle de jets méromorphe ${\sf J} \big/
y^c \, z^{ m(d-1)} \, t^{ m ( e - 1)}$ est aussi holomorphe au
voisinage de tout point $(x_0, y_0, z_0, 0)$ appartenant à $\{ z \neq
0\} \cap \{ t = 0\}$.

\medskip

Il reste à étudier (Cas 3) l'holomorphie de ${\sf J} \big/ y^c \, z^{
m(d-1)} \, t^{ m ( e - 1)}$ au voisinage d'un point appartenant à
l'ensemble fini des $d\, e$ points:
\[
\{z=0\}\cap\{t=0\}\cap X^2,
\]
qui s'identifient aux $d\, e$ points d'intersection entre
les deux courbes de $\C_{xy}^2$:
\[
\big\{R(x,y)=0\big\}
\ \ \ \ \ \ \ \ \ \ \ \ \ \ \ \ \
\text{\rm et}
\ \ \ \ \ \ \ \ \ \ \ \ \ \ \ \ \
\big\{S(x,y)=0\big\}.
\] 
Comme on s'y attend, cela va être garanti par une nouvelle application
du théorème de Riemann-Hurwitz-Hartogs, puisque ces points sont encore
de codimension $2$ dans la surface.

\medskip\noindent{\bf Disposition géométrique générique.}
{\em L'intersection:}
\[
X^2\cap\big\{z=t=0\big\},
\]
{\em qui s'identifie à l'intersection entre les deux courbes
de $\C_{ x, y}^2$:}
\[
\big\{
R
=
0
\big\}
\cap
\big\{
S
=
0
\big\},
\]
{\em ensemble fini de cardinal $d \cdot e$, 
ne contient aucun point de l'axe $Ox$.}

\medskip
Alors en un tel point: 
\[
(x_0,y_0,0,0)
\,\in\,
X^2\cap\big\{z=t=0\big\},
\]
pour lequel $y_0 \neq 0$ par disposition géométrique générique, 
la matrice jacobienne a ses deux dernières colonnes nulles:
\[
\left(\!\!
\begin{array}{cccc}
-R_x^0 & -R_y^0 & 0 & 0
\\
-S_x^0 & -S_y^0 & 0 & 0
\end{array}
\!\!\right),
\]
et comme les courbes sont à croisement normal en tout
point où elles s'intersectent, le déterminant
des deux premières colonnes est nécessairement non nul en ce point.
Le théorème des fonctions implicites assure alors
l'existence de deux fonctions holomorphes ${\tt X}$ et
${\tt Y}$ définies au voisinage de $(0, 0) = (z_0, t_0)$ dans 
$\C^2$ satisfaisant:
\[
\aligned
{\tt X}(0,0)
&
=
x_0,
\\
{\tt Y}(0,0)
&
=
y_0
\neq
0
\endaligned
\]
qui représentent $X^2$ localement comme le graphe:
\[
\aligned
x
&
=
{\tt X}(z,t),
\\
y
&
=
{\tt Y}(z,t).
\endaligned
\]

Maintenant, dans le quotient:
\[
\!\!\!\!\!\!\!\!\!\!\!\!\!\!\!\!\!\!\!\!
\frac{{\sf J}}{y^c\,z^{m(d-1)}\,t^{m(e-1)}}
=
\frac{1}{y^c\,z^{m(d-1)}\,t^{m(e-1)}}\,
\sum_{j+k+p+q=m}\,
A_{j,k,p,q}\,
(x')^j\,(y')^k\,
\big(R'\big)^p\,\big(S'\big)^q\,
\big(R\big)^{m-p}\,\big(S\big)^{m-q},
\]
remplaçons alternativement $R$, $R'$ puis $S$, $S'$ par les valeurs
tirées de:
\[
\left[
\aligned
z^d
&
=
R,
\\
z'\,d\,z^{d-1}
&
=
R',
\endaligned\right.
\ \ \ \ \ \ \ \ \ \ \ \ \ \ \ \ \ \ \ \ \ \ \ \ \ \ \ \ \ \ \ \ \ \ \ \ 
\ \ \ \
\left[
\aligned
t^e
&
=
S,
\\
t'\,e\,t^{e-1}
&
=
S',
\endaligned\right.
\]
ce qui donne deux expressions du quotient, la première ayant
fait disparaître la singularité en $z$:
\[
\aligned
\frac{{\sf J}}{y^c\,z^{m(d-1)}\,t^{m(e-1)}}
=
\frac{1}{y^c\,t^{m(e-1)}}\,
\sum_{j+k+p+q=m}\,
&
d^p\,
A_{j,k,p,q}\,
\\
&
z^{m-p}\,(z')^p
\\
&
(x')^j\,(y')^k
\\
&
\big(S'\big)^q\,\big(S\big)^{m-q},
\endaligned
\]
la seconde ayant fait disparaître la singularité en $t$:
\[
\aligned
\frac{{\sf J}}{y^c\,z^{m(d-1)}\,t^{m(e-1)}}
=
\frac{1}{y^c\,z^{m(d-1)}}\,
\sum_{j+k+p+q=m}\,
&
e^q\,
A_{j,k,p,q}\,
\\
&
t^{m-q}\,(t')^q
\\
&
(x')^j\,(y')^k
\\
&
\big(R'\big)^p\,\big(R\big)^{m-p}.
\endaligned
\]
L'écriture intrinsèque à $X^2$ exige de remplacer toutes
les occurences de $x, y, x', y'$ par leurs valeurs
déduites des équations graphantes, et on obtient 
comme première expression complète:
\[
\!\!\!\!\!\!\!\!\!\!\!\!\!\!\!\!\!\!\!\!
\aligned
&
\frac{{\sf J}}{y^c\,z^{m(d-1)}\,t^{m(e-1)}}
\bigg\vert_{X^2\,\,\text{\rm près de}\,\,
(x_0,y_0,0,0)}
=
\\
&
=
\frac{1}{
\underbrace{{\tt Y}(z,t)^c}_{
\text{\sf holomorphe, car}
\atop
{\tt Y}(0,0)\,=\,y_0\,\neq\,0}
\,t^{m(e-1)}}\,
\sum_{j+k+p+q=m}\,
d^p\,
A_{j,k,p,q}
\big({\tt X}(z,t),{\tt Y}(z,t)\big)
\\
&
z^{m-p}\,(z')^p
\\
&
\big(z'\,{\tt X}_z(z,t)+t'\,{\tt X}_t(z,t)\big)^j
\\
&
\big(z'\,{\tt Y}_z(z,t)+t'\,{\tt Y}_t(z,t)\big)^k
\\
&
\Big(
\big(z'\,{\tt X}_z(z,t)+t'\,{\tt X}_t(z,t)\big)\,
S_x\big({\tt X}(z,t),{\tt Y}(z,t)\big)
+
\big(z'\,{\tt Y}_z(z,t)+t'\,{\tt Y}_t(z,t)\big)\,
S_y\big({\tt X}(z,t),{\tt Y}(z,t)\big)
\Big)^q
\\
&
\Big(S\big({\tt X}(z,t),{\tt Y}(z,t)\big)\Big)^{m-q},
\endaligned
\]
et comme deuxième expression complète:
\[
\!\!\!\!\!\!\!\!\!\!\!\!\!\!\!\!\!\!\!\!
\aligned
&
\frac{{\sf J}}{y^c\,z^{m(d-1)}\,t^{m(e-1)}}
\bigg\vert_{X^2\,\,\text{\rm près de}\,\,
(x_0,y_0,0,0)}
=
\\
&
=
\frac{1}{
\underbrace{{\tt Y}(z,t)^c}_{
\text{\sf holomorphe, car}
\atop
{\tt Y}(0,0)\,=\,y_0\,\neq\,0}
\,z^{m(d-1)}}\,
\sum_{j+k+p+q=m}\,
e^q\,
A_{j,k,p,q}
\big({\tt X}(z,t),{\tt Y}(z,t)\big)
\\
&
t^{m-q}\,(t')^q
\\
&
\big(z'\,{\tt X}_z(z,t)+t'\,{\tt X}_t(z,t)\big)^j
\\
&
\big(z'\,{\tt Y}_z(z,t)+t'\,{\tt Y}_t(z,t)\big)^k
\\
&
\Big(
\big(z'\,{\tt X}_z(z,t)+t'\,{\tt X}_t(z,t)\big)\,
R_x\big({\tt X}(z,t),{\tt Y}(z,t)\big)
+
\big(z'\,{\tt Y}_z(z,t)+t'\,{\tt Y}_t(z,t)\big)\,
R_y\big({\tt X}(z,t),{\tt Y}(z,t)\big)
\Big)^p
\\
&
\Big(R\big({\tt X}(z,t),{\tt Y}(z,t)\big)\Big)^{m-p}.
\endaligned
\]

Analysons ces deux résultats. Dans les coordonnées intrinsèques $(z,
t)$ sur $X^2$ près de $(x_0, y_0, 0, 0)$, la première expression est
singulière précisément sur $\{ t = 0\}$, et elle est holomorphe
ailleurs.  La seconde est singulière précisément sur $\{ z = 0\}$, et
elle est holomorphe ailleurs.  Le théorème d'élimination des
singularités de codimension $2$ 
de Riemann-Hurwitz-Hartogs
s'applique à nouveau dans cette
circonstance, et montre que la différentielle de jets méromorphe ${\sf
J} \big/ y^c\, z^{ m( d-1)} \, t^{ m ( e-1)}$ est aussi holomorphe
en chacun des points d'intersection de $X^2 \cap \P_{xy}^2$.  

\smallskip

La démonstration détaillée de la première
proposition principale est donc achevée.
\endproof

\noindent{\bf Résumé.} On vient d'établir que 
la restriction à la partie affine $X^2 \cap \C^4$ de la
différentielle de jets méromorphe extrinsèque
${\sf J} \big/ 
y^c\, z^{ m( d-1)} \, t^{ m ( e-1)}$ est holomorphe, et donc,
il ne reste plus qu'à étudier l'holomorphie au voisinage
de tout point appartenant à $X^2 \cap \P_\infty^3$.

\medskip\noindent{\bf Lemme.}
{\em
L'intersection de $X^2$ avec la droite à l'infini $\P_{ \infty, zt}^1$
dans la direction des deux axes $0z$ et $0t$ est réduite à l'ensemble
vide; de plus, le $(1/x)$-
et le $(1 / y)$-changements de carte affine
suffisent à capturer tous les points de la courbe à l'infini
$X^2 \cap \P_{\infty, xyzt}^3$.}

\proof
En effet, rappelons qu'en termes des coordonnées homogènes 
$[U \colon X \colon Y \colon Z \colon T]$ sur $\P^4$, 
les coordonnées affines initiales sont égales à:
\[
\big[
1
\colon
{\textstyle{\frac{X}{U}}},\,
{\textstyle{\frac{Y}{U}}},\,
{\textstyle{\frac{Z}{U}}},\,
{\textstyle{\frac{T}{U}}}
\big]
\equiv
\big(1\colon x\colon y\colon z\colon t\big).
\]
Elles couvrent l'ouvert $\{ U \neq 0\}$.
Le $(1/x)$-changement de carte:
\[
(x,y,z,t)
\,\longmapsto\,
\big(
{\textstyle{\frac{1}{x}}},\,
{\textstyle{\frac{y}{x}}},\,
{\textstyle{\frac{z}{x}}},\,
{\textstyle{\frac{t}{x}}}
\big)
=:
(x_1,y_1,z_1,t_1)
\]
couvre l'ouvert $\{ X \neq 0\}$, et le 
Le $(1/y)$-changement de carte:
\[
(x,y,z,t)
\,\longmapsto\,
\big(
{\textstyle{\frac{x}{y}}},\,
{\textstyle{\frac{1}{y}}},\,
{\textstyle{\frac{z}{y}}},\,
{\textstyle{\frac{t}{y}}}
\big)
=:
(x_2,y_2,z_2,t_2)
\]
couvre l'ouvert $\{ Y \neq 0\}$. Seuls les points de la droite: 
\[
\big\{
[
0
\colon
0
\colon
0
\colon
Z
\colon
T
]
\colon\,
[Z\colon T]\,\in\,\P^1
\big\}
\]
contenue dans le plan $\P_{\infty,xyzt}^3$ à l'infini ne sont pas
couverts par ces trois cartes affines.  Or nous affirmons qu'aucun de
ces points ne peut appartenir à notre surface $X^2 \subset \P^4$.

En effet, les deux équations homogènes de $X^2$
s'obtiennent simplement en
homogénéisant ses deux équations affines:
\[
\aligned
{\rm R}(U\colon X\colon Y\colon Z\colon T)
&
:=
Z^d
-
U^d\,R\big(
{\textstyle{\frac{X}{U}}},\,
{\textstyle{\frac{X}{U}}}
\big)
\equiv
Z^d
-
{\rm R}^\sim(U\colon X\colon Y),
\\
{\rm S}(U\colon X\colon Y\colon Z\colon T)
&
:=
T^e
-
U^d\,S\big(
{\textstyle{\frac{X}{U}}},\,
{\textstyle{\frac{X}{U}}}
\big)
\equiv
T^e
-
{\rm S}^\sim(U\colon X\colon Y).
\endaligned
\]
Les deux polynômes homogènes de trois variables ${\rm R}^\sim$ et
$S^\sim$ ainsi obtenus s'annulant automatiquement en $[0 \colon 0
\colon 0]$, on voit bien que $[0 \colon 0 \colon 0 \colon Z \colon T
] \in X^2$ implique $Z^d = 0$ et $T^e = 0$, contredisant le fait que
$[Z \colon T] \in \P^1$ ne doit jamais avoir ses deux coordonnées
homogènes nulles.
\endproof

Dans l'espace projectif: 
\[
\P^4(\C)
=
\big\{
[U\colon X\colon Y\colon Z\colon T]
\big\}, 
\]
si l'on part des coordonnées affines {\sl initiales}:
\[
x
:=
\frac{X}{U},
\ \ \ \ \ \ \ \ \ \ \ \ \ \ \ \ \
y
:=
\frac{Y}{U},
\ \ \ \ \ \ \ \ \ \ \ \ \ \ \ \ \
z
:=
\frac{Z}{U},
\ \ \ \ \ \ \ \ \ \ \ \ \ \ \ \ \
t
:=
\frac{T}{U},
\]
les quatre autres sytèmes de coordonnées affines
s'en déduisent par:

\medskip\noindent$\bullet$\,
le $(1/x)$-changement de carte:
\[
\aligned
\big(x,y,z,t\big)
&
\,\longmapsto\,
\bigg(
\frac{1}{x},\,
\frac{y}{x},\,
\frac{z}{x},\,
\frac{t}{x}
\bigg)
\\
&\ \ \ \,
=:
\big(x_1,y_1,z_1,t_1\big);
\endaligned
\]

\medskip\noindent$\bullet$\,
le $(1/y)$-changement de carte:
\[
\aligned
\big(x,y,z,t\big)
&
\,\longmapsto\,
\bigg(
\frac{x}{y},\,
\frac{1}{y},\,
\frac{z}{y},\,
\frac{t}{y}
\bigg)
\\
&\ \ \ \,
=:
\big(x_2,y_2,z_2,t_2\big);
\endaligned
\]

\medskip\noindent$\bullet$\,
le $(1/z)$-changement de carte:
\[
\aligned
\big(x,y,z,t\big)
&
\,\longmapsto\,
\bigg(
\frac{x}{z},\,
\frac{y}{z},\,
\frac{1}{z},\,
\frac{t}{z}
\bigg)
\\
&\ \ \ \,
=:
\big(x_3,y_3,z_3,t_3\big);
\endaligned
\]

\medskip\noindent$\bullet$\,
le $(1/t)$-changement de carte:
\[
\aligned
\big(x,y,z,t\big)
&
\,\longmapsto\,
\bigg(
\frac{x}{t},\,
\frac{y}{t},\,
\frac{z}{t},\,
\frac{1}{t}
\bigg)
\\
&\ \ \ \,
=:
\big(x_4,y_4,z_4,t_4\big).
\endaligned
\]

En partant de la carte affine initiale
$(x, y, z, t) \in \C^4$, ce qu'on couvre 
avec le $(1/x)$-changement de carte
et le $(1/y)$-changement de carte, c'est justement:
\[
\P^4
\big\backslash
\P_{\infty,zt}^1.
\]

{\em Par conséquent, pour analyser le comportement
à l'infini de la différentielle de jets
extrinsèque méromorphe, il suffit d'examiner
comment elle se transfère à travers ces
deux seuls changements de carte affine.}

Par symétrie formelle entre les coordonnées $(x, y)$,
il suffit même d'ailleurs d'examiner le $(1/x)$-changement
de carte.

\`A travers, donc:
\[
\aligned
\big(x,y,z,t\big)
&
\,\longmapsto\,
\bigg(
\frac{1}{x},\,
\frac{y}{x},\,
\frac{z}{x},\,
\frac{t}{x}
\bigg)
\\
&\ \ \ \,
=:
\big(x_1,y_1,z_1,t_1\big),
\endaligned
\]
tout d'abord les deux polynômes $R$ et $S$ sont naturellement
transformés
en deux nouveaux polynômes $R_1$ et $S_1$ définis par:
\[
\aligned
R_1(x_1,y_1)
&
:=
(x_1)^d\,R
\bigg(
\frac{1}{x_1},\,\frac{y_1}{x_1}
\bigg),
\\
S_1(x_1,y_1)
&
:=
(x_1)^e\,S
\bigg(
\frac{1}{x_1},\,\frac{y_1}{x_1}
\bigg).
\endaligned
\]

Pour calculer comment se transforme:
\[
\!\!\!\!\!\!\!\!\!\!\!\!\!\!\!\!\!\!\!\!
\aligned
J\big(x,y,x',y'\big)
=
\sum_{j+k+p+q=m}\,
A_{j,k,p,q}(x,y)\,
\big(x'\big)^j\,
\big(y'\big)^k\,
&
\Big(
x'\,R_x(x,y)
+
y'\,R_y(x,y)
\Big)^p\,
\big(R(x,y)\big)^{m-p}
\\
&
\Big(
x'\,S_x(x,y)
+
y'\,S_y(x,y)
\Big)^q\,
\big(S(x,y)\big)^{m-q},
\endaligned
\]
le transfert de nombreux autres termes doit être examiné
séparément à l'avance.

Tout d'abord:
\[
\aligned
x'
&
=
\bigg(
\frac{1}{x_1}
\bigg)'
=
-\,\frac{x_1'}{x_1x_1},
\\
y'
&
=
\bigg(
\frac{y_1}{x_1}
\bigg)'
=
\frac{y_1'}{x_1}
-
\frac{x_1'y_1}{x_1x_1}.
\endaligned
\]

Ensuite, si on différentie:
\[
R(x,y)
=
(x)^d\,
R_1
\bigg(
\frac{1}{x},\,\frac{y}{x}
\bigg),
\]
premièrement par rapport à $x$, on obtient:
\[
R_x(x,y)
=
d\,(x)^{d-1}\,
R_1
\bigg(
\frac{1}{x},\,\frac{y}{x}
\bigg)
-
\frac{(x)^d}{xx}\,
R_{1,x_1}
\bigg(
\frac{1}{x},\,\frac{y}{x}
\bigg)
-
\frac{(x)^dy}{xx}\,
R_{1,y_1}
\bigg(
\frac{1}{x},\,\frac{y}{x}
\bigg),
\]
et deuxièmement par rapport à $y$, on obtient:
\[
R_y(x,y)
=
\frac{(x)^d}{x}\,
R_{1,y_1}
\bigg(
\frac{1}{x},\,\frac{y}{x}
\bigg).
\]

Tout cela mis ensemble donne, après simplifications et regroupements:
\[
\!\!\!\!\!\!\!\!\!\!\!\!\!\!\!\!\!\!\!\!
\aligned
x'R_x(x,y)
+
y'R_y(x,y)
&
=
\frac{1}{(x_1)^{d+1}}\,
\Big[
-\,x_1'\,d\,R_1(x_1,y_1)
\Big]
+
\frac{1}{(x_1)^d}\,
\Big[
x_1'R_{1,x_1}(x_1,y_1)
+
y_1'R_{1,y_1}(x_1,y_1)
\Big]
\\
&
=
-\,
\frac{x_1'}{(x_1)^{d+1}}\,d\,
R_1(x_1,y_1)
+
\frac{1}{(x_1)^d}\,
R_1',
\endaligned
\]
équation que l'on peut d'ailleurs aussi trouver directement en
différentiant (prendre le prime) formellement:
\[
R
=
\frac{1}{(x_1)^d}\,R_1.
\]

Partant de:
\[
S
=
\frac{1}{(x_1)^e}\,S_1,
\]
on trouve de manière entièrement similaire:
\[
S'
=
-\,\frac{x_1'}{(x_1)^{e+1}}\,e\,S_1
+
\frac{1}{(x_1)^e}\,S_1'.
\]

Tous ces calculs préparatoires élémentaires nous permettent d'écrire
maintenant le transfert brut de la différentielle de jets méromorphe
extrinsèque à travers le $(1/x)$-changement de carte:
\[
\!\!\!\!\!\!\!\!\!\!\!\!\!\!\!\!\!\!\!\!
\aligned
\frac{J(x,y,x',y')}{
y^c\,z^{m(d-1)}\,t^{m(e-1)}}
&
=
\frac{\big(x_1\big)^{c+m(d-1)+m(e-1)}}{
(y_1)^c\,(z_1)^{m(d-1)}\,(t_1)^{m(e-1)}}\,
\sum_{j+k+p+q=m}\,
\sum_{h+i\leqslant a}\,
\overbrace{A_{j,k,p,q}^{h,i}}^{{\sf coefficients}\,\in\,\C}
\frac{1}{(x_1)^h}\,
\bigg(
\frac{y_1}{x_1}
\bigg)^i
\\
&
\ \ \ \ \
\bigg(\!
-\,
\frac{x_1'}{x_1x_1}
\bigg)^j\,
\bigg(
\frac{x_1y_1'-y_1x_1'}{x_1x_1}
\bigg)^k\,
\\
&
\ \ \ \ \
\bigg(\!
-\,\frac{x_1'\,d}{(x_1)^{d+1}}\,
R_1(x_1,y_1)
+
\frac{1}{(x_1)^d}\,
\Big[
x_1'\,R_{1,x_1}(x_1,y_1)
+
y_1'\,R_{1,y_1}(x_1,y_1)
\Big]
\bigg)^p
\\
&
\ \ \ \ \
\bigg(\!
-\,\frac{x_1'\,e}{(x_1)^{e+1}}\,
S_1(x_1,y_1)
+
\frac{1}{(x_1)^e}\,
\Big[
x_1'\,S_{1,x_1}(x_1,y_1)
+
y_1'\,S_{1,y_1}(x_1,y_1)
\Big]
\bigg)^q
\\
&
\ \ \ \ \
\bigg(
\frac{1}{(x_1)^d}\,R_1(x_1,y_1)
\bigg)^{m-p}\,
\bigg(
\frac{1}{(x_1)^e}\,S_1(x_1,y_1)
\bigg)^{m-q}.
\endaligned
\]
Extrayons alors la présence de $x_1$ dans les nombreux dénominateurs
comme suit:
\[
\aligned
\frac{J(x,y,x',y')}{
y^c\,z^{m(d-1)}\,t^{m(e-1)}}
&
=
\frac{\big(x_1\big)^{c+md-m+me-d}}{
(y_1)^c\,(z_1)^{m(d-1)}\,(t_1)^{m(e-1)}}\,
\sum_{j+k+p+q=m}\,
\sum_{h+i\leqslant a}\,
A_{j,k,p,q}^{h,i}\,
\big(y_1\big)^i
\\
&
\ \ \ \ \
\frac{1}{\big(x_1\big)^{
h+i+2j+2k+(d+1)p+(e+1)q+d(m-p)+e(m-q)}}\,
\big(\!-x_1'\big)^j\,
\big(x_1y_1'-y_1x_1'\big)^k
\\
&
\ \ \ \ \
\Big(\!
-\,x_1'\,d\,R_1(x_1,y_1)
+
x_1
\big[
x_1'\,R_{1,x_1}(x_1,y_1)
+
y_1'\,R_{1,y_1}(x_1,y_1)
\big]
\Big)^p
\\
&
\ \ \ \ \
\Big(\!
-\,x_1'\,e\,S_1(x_1,y_1)
+
x_1
\big[
x_1'\,S_{1,x_1}(x_1,y_1)
+
y_1'\,S_{1,y_1}(x_1,y_1)
\big]
\Big)^q
\\
&
\ \ \ \ \
\Big(
R_1(x_1,y_1)
\Big)^{m-p}\,
\Big(
S_1(x_1,y_1)
\Big)^{m-q}.
\endaligned
\]
Pour simplifier encore la présence de $x_1$,
servons-nous du $x_1$ situé au numérateur 
que l'on décompose à cette effet sous la forme:
\[
\big(x_1\big)^{c+md-m+em-m}
=
\big(x_1\big)^{c-a-4m}\,
\underbrace{\big(x_1\big)^{a+dm+m+me+m}}_{
\text{\sf à placer après $\sum\,\sum$}},
\]
ce qui conduit à:
\[
\frac{\big(x_1\big)^{a+dm+m+me+m}}{
\big(x_1\big)^{h+i+2j+2k+dp+p+eq+q+dm-dp+em-eq}}
=
\big(x_1\big)^{a-(h+i)+2m-(2j+2k+p+q)},
\]
et cette dernière puissance de $x_1$ est automatiquement
$\geqslant 0$, puisque par hypothèse:
\[
\aligned
h+i
&
\leqslant
a,
\\
2j+2k+p+q
&
\leqslant
2j+2k+2p+2q
\\
&
=
2m.
\endaligned
\]

En définitive, le transfert à l'infini devient:
\[
\aligned
\frac{J(x,y,x',y')}{
y^c\,z^{m(d-1)}\,t^{m(e-1)}}
&
=
\frac{\big(x_1\big)^{c-a-4m}}{
(y_1)^c\,(z_1)^{m(d-1)}\,(t_1)^{m(e-1)}}\,
\sum_{j+k+p+q=m}\,
\sum_{h+i\leqslant a}\,
A_{j,k,p,q}^{h,i}\,
\big(y_1\big)^i
\\
&
\ \ \ \ \
\big(x_1\big)^{a-(h+i)+2m-(2j+2k+p+q)}\,
\big(\!-x_1'\big)^j\,
\big(x_1y_1'-y_1x_1'\big)^k
\\
&
\ \ \ \ \
\Big(\!
-\,x_1'\,d\,R_1(x_1,y_1)
+
x_1
\big[
x_1'\,R_{1,x_1}(x_1,y_1)
+
y_1'\,R_{1,y_1}(x_1,y_1)
\big]
\Big)^p
\\
&
\ \ \ \ \
\Big(\!
-\,x_1'\,e\,S_1(x_1,y_1)
+
x_1
\big[
x_1'\,S_{1,x_1}(x_1,y_1)
+
y_1'\,S_{1,y_1}(x_1,y_1)
\big]
\Big)^q
\\
&
\ \ \ \ \
\Big(
R_1(x_1,y_1)
\Big)^{m-p}\,
\Big(
S_1(x_1,y_1)
\Big)^{m-q}.
\endaligned
\]

La partie affine $\C_{\infty,x}^3$ à l'infini dans
la direction de l'axe $Ox$ est maintenant devenue:
\[
\big\{
x_1
=
0
\big\}.
\]
\'Etudier l'holomorphie à l'infini revient donc
à étudier l'holomorphie sur ce $3$-plan $\{ x_1 = 0 \}$.
(L'étude de ce qui se passe à travers le
$(1/y)$-changement de carte est entièrement similaire.)

Observons que tous les termes après la double
$\sum\, \sum$ sont maintenant devenus
{\em polynomiaux} en $(x_1, y_1, x_1', y_1')$,
donc holomorphes.
Seule la fraction qui les précède:
\[
\frac{\big(x_1\big)^{c-a-4m}}{
(y_1)^c\,(z_1)^{m(d-1)}\,(t_1)^{m(e-1)}}
\]
peut présenter des pôles.

Or maintenant, nous allons supposer expressément que:
\[
c-a-4m
\,\geqslant\,0,
\]
ce qui revient à restreindre la borne $a$ sur le degré 
des $A_{ j, k, p, q}$ à satisfaire:
\[
\boxed{\,
a
\,\leqslant\,
c
-
4m.\,}
\]
\`A un décalage d'une unité près, si on suppose même que:
\[
a
\,\leqslant\,
c
-
4m-1,
\]
la totalité de l'expression est multiple de $x_1$,
et donc, elle s'annule sur $\{ x_1 = 0\}$ en dehors du lieu polaire.

\medskip\noindent{\bf Diposition géométrique générique.}
{\em On peut supposer à l'avance que:}

\smallskip\noindent$\bullet$\,
{\em l'intersection:}
\[
X^2
\cap
\P_{\infty,xyzt}^3
\]
{\em est une courbe lisse, dont la partie
$X^2 \cap \C_{\infty, x}^3$ est devenue ici
$\{ x_1 = 0\} \cap X^2$};

\smallskip\noindent$\bullet$\,
{\em dans ce $3$-plan complexe, les trois intersections
suivantes de $X^2$
avec des $2$-plans complexes:}
\[
\aligned
&
X^2\cap\{0=x_1=y_1\},
\\
&
X^2\cap\{0=x_1=z_1\},
\\
&
X^2\cap\{0=x_1=t_1\},
\endaligned
\]
{\em se réduisent chacune à un nombre fini de points distincts
deux à deux.}

\medskip

L'expression obtenue à l'instant s'avère alors être holomorphe
en tout point $\big(0, y_{10}, z_{10}, t_{10} \big)$
appartenant à l'ensemble:
\[
\aligned
X^2
&
\cap
\big\{x_1=0\big\}
\\
&
\cap
\big\{y_1\neq 0\big\}
\cap
\big\{z_1\neq 0\big\}
\cap
\big\{t_1\neq 0\big\}.
\endaligned
\]
Puisque l'ensemble restant: 
\[
\Big(
X^2\cap\big\{0=x_1=y_1\big\}
\Big)
\bigcup
\Big(
X^2\cap\big\{0=x_1=z_1\big\}
\Big)
\bigcup
\Big(
X^2\cap\big\{0=x_1=t_1\big\}
\Big)
\]
est de codimension $2$
dans $X^2$, le théorème d'élimination des singularités
de Riemann-Hurwitz-Hartogs conclut la 
démonstration de l'énoncé suivant.

\medskip\noindent{\bf Deuxième proposition principale.}
{\em
Sur une surface $X^2 \subset \P^4$ d'équations affines:}
\[
\aligned
z^d
&
=
R(x,y),
\\
t^e
&
=
S(x,y),
\endaligned
\]
{\em où $R \in \C [x, y]$ et $S \in \C [x, y]$ sont de degrés
respectifs $d \leqslant e$ qui satisfont un
nombre fini de dispositions géométriques
génériques,
pour tous polynômes:}
\[
A_{j,k,p,q}(x,y)
\,\in\,\C[x,y]
\]
{\em de degrés:}
\[
\deg\,
A_{j,k,p,q}
\,\leqslant\,
a
\leqslant
c-4m,
\]
{\em la différentielle de jets extrinsèque méromorphe:}
\[
\aligned
\frac{\sf J}{y^c\,z^{m(d-1)}\,t^{m(e-1)}}
&
=
\frac{{\sf J}(x,y,x',y')}{y^c\,z^{m(d-1)}\,t^{m(e-1)}}
\\
&
=
\frac{1}{y^c\,z^{m(d-1)}\,t^{m(e-1)}}\,
\sum_{j+k+p+q=m}\,
A_{j,k,p,q}\,
\big(x'\big)^j\,\big(y'\big)^k\,
\\
&
\ \ \ \ \ \ \ \ \ \ \ \ \ \ \ \ \ \ \ \ \ \ \ \ \ \ \ \ \ \ \ \ \ \ \ \ 
\ \ \ \ \ \ \ \ \ \ \ \ \ \ \
\Big(
x'\,R_x+y'\,R_y
\Big)^p\,
\Big(
x'\,S_x+y'\,S_y
\Big)^q
\\
&
\ \ \ \ \ \ \ \ \ \ \ \ \ \ \ \ \ \ \ \ \ \ \ \ \ \ \ \ \ \ \ \ \ \ \ \ 
\ \ \ \ \ \ \ \ \ \ \ \ \ \ \
\big(R\big)^{m-p}\,
\big(S\big)^{m-q}
\endaligned
\]
{\em possède une restriction à $X^2$:}
\[
\frac{\sf J}{y^c\,z^{m(d-1)}\,t^{m(e-1)}}
\bigg\vert_{X^2}
\]
{\em qui est une section {\em holomorphe} du fibré des différentielles
symétriques intrinsèques:}
\[
\Sym^m T_X^*,
\]
{\em pourvu seulement que:}
\[
{\sf J}\big(x,y,x',y'\big)
\,\,\equiv\,\,
y^c\,\,\widetilde{\sf J}\big(x,y,x',y'\big)
\]
{\em soit divisible par $y^c$.}

\medskip

C'est à l'analyse de cette contrainte de divisibilité
qu'est consacrée la dernière Section~5 qui suit.


\bigskip

\centerline{\bf 5.~Contraintes et degrés de liberté}
\label{contraintes-liberte}
\HEAD{5.~Contraintes et degrés de liberté}{
Jo\"el Merker, Département de Mathématiques d'Orsay}

\medskip

Pour fixer les idées, prenons maintenant:
\[
\boxed{\,
c
:=
d,\,}
\]
ainsi que:
\[
\boxed{\,
a
:=
d-4\,m.}
\]

Chaque polynôme:
\[
A_{j,k,p,q}
=
\sum_{h+i\leqslant a}\,
A_{j,k,p,q}^{h,i}\,
x^h\,y^i
\]
de degré $\leqslant a$ possède:
\[
\frac{(a+1)(a+2)}{1\cdot 2}
\]
coefficients complexes libres $A_{j,k,p,q}^{h,i}$.
Puisque:
\[
\frac{(m+1)(m+2)(m+3)}{1\cdot 2\cdot 3}
=
\Card\,
\Big\{
(j,k,p,q)\in\N^4\colon\,
j+k+p+q=m
\Big\},
\]
il y a au total:
\[
\frac{(a+1)(a+2)}{1\cdot 2}\,
\frac{(m+1)(m+2)(m+3)}{1\cdot 2\cdot 3}
\]
degrés de liberté dans les $A_\smallbullet$, et nous retiendrons
la minoration:
\[
\boxed{\,
\mathmotsf{degrés de liberté dans les}\,\,
A_\smallbullet
\,\geqslant\,
\frac{(d-4m)^2}{2}\,
\frac{m^3}{6}.\,}
\]

Lorsqu'on développe et réorganise le numérateur ${\sf J}$, on
obtient comme cela a déjà été vu p.~\pageref{deja-vu}:
\[
{\sf J}
=
\sum_{\alpha+\beta\,=\,m}\,
(x')^\alpha\,(y')^\beta\,
\Lambda_{\alpha,\beta}
\big(A_\smallbullet;\,
R_x,R_y,S_x,S_y,R,S
\big),
\]
où chaque coefficient:
\[
\aligned
\Lambda_{\alpha,\beta}
&
=
\sum_{j+p_1+q_1\,=\,\alpha}\,\,\,
\sum_{k+p_2+q_2\,=\,\beta}\,
{\textstyle{\frac{(p_1+p_2)!}{p_1!\,p_2!}}}\,
{\textstyle{\frac{(q_1+q_2)!}{q_1!\,q_2!}}}\,
A_{j,k,p_1+p_2,q_1+q_2}\,
\\
&
\ \ \ \ \ \ \ \ \ \ \ \ \ \ \ \ \ \ \ \ \ \ \ \ \ \ \ \ \ \ \ \ \ \ 
\ \ \ \ \ \ \ \ \ \ 
\big(R_x\big)^{p_1}\,\big(R_y\big)^{p_2}\,
\big(S_x\big)^{q_1}\,\big(S_y\big)^{q_2}\,
\\
&
\ \ \ \ \ \ \ \ \ \ \ \ \ \ \ \ \ \ \ \ \ \ \ \ \ \ \ \ \ \ \ \ \ \
\ \ \ \ \ \ \ \ \ \  
\big(R\big)^{m-p_1-p_2}\,\big(S\big)^{m-q_1-q_2},
\endaligned
\]
visiblement linéaire par rapport aux $A_\smallbullet$, 
est un polynôme en $(x, y)$ de degré (exercice mental):
\[
\aligned
\deg\,\Lambda_{\alpha,\beta}
&
\,\leqslant\,
a
+
(d-1)\,p
+
d\,(m-p)
+
(e-1)\,q
+
e\,(m-q)
\\
&
\,=\,
a
+
dm-p
+
em-q
\\
&
\,\leqslant\,
d-1
+
dm
+
em.
\endaligned
\]
La divisibilité par $y^d$ d'un polynôme $\Lambda_{ \alpha, \beta}$
demande au plus l'annulation des coefficients des monômes:
\[
\aligned
&
x^0y^{d-1},\ \ x^1y^{d-1},\ \
\dots\dots,\ \
x^{dm+em}y^{d-1}
\\
&
\cdots\cdots\cdots
\\
&
x^0y^0,\ \ \ x^1y^0,\ \ \
\dots\dots\dots,\ \
x^{dm+em}y^0,\ \
\dots\dots,\ \
x^{d-1+dm+em}y^0,\ \
\endaligned
\]
donc le nombre d'équations linéaires (contraintes)
auxquelles sont soumis les $A_\bullet$ pour chaque $\Lambda_{ \alpha,
\beta}$ est majoré par (substituer au trapèze
d'exposants le rectangle qui lui est exinscrit):
\[
\leqslant\,
d\,\big[d+dm+em\big].
\]
Comme le nombre de polynômes $\Lambda_{ \alpha, \beta}$ est égal à:
\[
m+1,
\]
on conclut que:
\[
\boxed{\,
\mathmotsf{nombre de contraintes linéaires en les}\,\,
A_\smallbullet
\,\leqslant\,
(m+1)\,d\,
\big[d+dm+em\big].\,}
\]

Prenons maintenant pour quasiment optimiser:
\[
m
:=
\Ent\,
\bigg(
\frac{d}{12}
\bigg),
\]
de telle sorte que:
\[
\frac{d}{12}-1
\,\leqslant\,
m
\,\leqslant\,
\frac{d}{12}.
\]

Les degrés de libertés sont alors minorés par:
\[
\aligned
\mathmotsf{degrés de liberté}
&
\,=\,
\frac{(d-4m)^2}{2}\,\frac{m^3}{6}
\\
&
\,\geqslant\,
\frac{1}{12}\,
\bigg(
d
-
\frac{4\,d}{12}
\bigg)^2\,
\bigg(
\frac{d}{12}-1
\bigg)^3
\\
&
\,=\,
\frac{1}{12}\,\frac{4}{9}\,
d^2\,
\bigg(
\frac{d}{12}-1
\bigg)^3,
\endaligned
\]
tandis que les contraintes sont majorées par:
\[
\aligned
\bigg(
\frac{d}{12}+1
\bigg)\,
d\,
\bigg[
d+d\,\frac{d}{12}+e\,\frac{d}{12}
\bigg]
&
\,\geqslant\,
\underbrace{(m+1)\,d\,
\big[d+dm+em\big]}_{\sf contraintes}.
\endaligned
\]

Puisque tout système linéaire dont le nombre d'inconnues excède celui
d'équations admet automatiquement des solutions,
il suffit donc que ce premier minorant auxiliaire
majore ce second majorant auxiliaire, à savoir il suffit que:
\[
\frac{1}{12}\,\frac{4}{9}\,
d^2\,
\bigg(
\frac{d}{12}-1
\bigg)^3
\,\geqslant\,
\bigg(
\frac{d}{12}+1
\bigg)\,
d\,
\bigg[
d
+
d\,\frac{d}{12}
+
e\,\frac{d}{12}
\bigg],
\]

En supposant maintenant comme 
cela a été stipulé dans l'énoncé du théorème que:
\[
\aligned
d
\,\leqslant\,
e
&
\,\leqslant\,
\frac{1}{2}\,\frac{1}{27}\,\frac{1}{12}\,d^2
\\
&
\,=\,
\frac{1}{648}\,d^2,
\endaligned
\]
on assure que le coefficient rationnel de la puissance maximale $d^3$
est positif dans l'inégalité suffisante (après suppression de $d^2$):
\[
\frac{1}{12}\,\frac{4}{9}\,
\zero{d^2}\,
\bigg(
\frac{d}{12}-1
\bigg)^3
\,\geqslant\,
\bigg(
\frac{d}{12}+1
\bigg)\,
\zero{d^2}\,
\bigg[
1
+
\frac{d}{12}
+
\frac{d^2}{2\cdot 27\cdot 12\cdot 12}
\bigg],
\]
inégalité qui devient après placement à gauche de tous les termes:
\[
\frac{1}{93312}\,d^3
-
\frac{61}{7776}\,d^2
-
\frac{17}{108}\,d
-
\frac{28}{27}
\,\geqslant\,
0,
\]
inégalité qui s'avère être satisfaite à partir de:
\[
d
\,\geqslant\,
752,
\]
ce qui conclut la
démonstration du théorème.\qed

\vfill\end{document}

%% file: macros.tex
\newcommand{\C}{\mathbb{C}}

\newcommand{\N}{\mathbb{N}}\renewcommand{\P}{\mathbb{P}}

\newcommand{\red}{\textcolor{red}}

\newcommand{\HEAD}[2]{%
\pagestyle{fancy}
\fancyhead[RO]{\tiny\sf\thepage}
\fancyhead[CO]{{\tiny\sf #1}}
\fancyhead[LE]{\tiny\sf\thepage}
\fancyhead[CE]{{\tiny\sf #2}}
\fancyfoot{}}

\theoremstyle{definition}

\newcommand{\Card}{\text{\scriptsize\sf Card}}

\renewcommand{\deg}{\text{\footnotesize\sf deg}}

\renewcommand{\dim}{\text{\footnotesize\sf dim}}

\newcommand{\Ent}{\text{\footnotesize\sf Ent}}

\newcommand{\explication}[1]{\text{\scriptsize\sf [#1]}}

\renewcommand{\lim}{\text{\footnotesize\sf lim}}

\newcommand{\mathmotsf}[1]{\text{\footnotesize\sf #1}}

\newcommand{\smallbullet}{{\scriptscriptstyle{\bullet}}}

\newcommand{\Sym}{\text{\scriptsize\sf Sym}}

\newcommand{\vf}{\vfill\end{document}}

\newcommand{\zero}[1]{\underline{#1}_{\red{\circ}}}
